\documentclass[11pt]{article}

\usepackage{fullpage}
\usepackage{secdot}
\usepackage{epsfig}
\usepackage{amsmath}
\usepackage{amssymb}
\usepackage{amsbsy}
\usepackage{amsthm}
\usepackage{graphicx}
\usepackage{geometry}
\usepackage{algorithm}
\usepackage{algorithmicx}
\usepackage{authblk}

\newtheorem{Definition}{Definition}
\newtheorem{Assumption}{Assumption}
\newtheorem{Lemma}{Lemma}
\newtheorem{Theorem}{Theorem}

\newtheorem{Proposition}{Proposition}

\newcommand{\iid}{\stackrel{\mathrm{iid}}{\sim}}
\newcommand{\parenthnewln}{\right.\\&\left.\quad\quad{}}
\usepackage{color}
\def\ez#1{{\color{black} #1}}

\begin{document}

\title{Gradient-Based Adaptive Stochastic Search for Non-Differentiable Optimization}
\date{First draft: October 21, 2011\\
This version: October 22, 2012}

\author{Enlu Zhou}
\affil{\scriptsize{Department of Industrial $\&$
Enterprise Systems Engineering,
University of Illinois at
Urbana-Champaign, IL 61801, enluzhou@illinois.edu}}
\author{Jiaqiao Hu}
\affil{\scriptsize{Department of Applied Mathematics and Statistics, Stony Brook University, NY 11794, jqhu@ams.sunysb.edu}}

\maketitle
\section*{ABSTRACT}
In this paper, we propose a stochastic search algorithm for solving general optimization problems with little structure. The algorithm iteratively finds high quality solutions by randomly sampling candidate solutions from a parameterized distribution model over the solution space. The basic idea is to convert the original (possibly non-differentiable) problem into a differentiable optimization problem on the parameter space of the parameterized sampling distribution, and then use a direct gradient search method to find improved sampling distributions. Thus, the algorithm combines the robustness feature of stochastic search from considering a population of candidate solutions with the relative fast convergence speed of classical gradient methods by exploiting local differentiable structures. We analyze the convergence and converge rate properties of the proposed algorithm, and carry out numerical study to illustrate its performance.

\section{Introduction}

We consider global optimization problems over real vector-valued domains. These optimization problems arise in many areas of importance and can be extremely difficult to solve due to
the presence of multiple local optimal solutions and the lack of structural properties such as differentiability and convexity. In such a general setting, there is little problem-specific knowledge that can be exploited in searching for improved solutions, and it is often the case that the objective function can only be assessed through the form of ``black-box" evaluation, which returns the function value for a specified candidate solution.

An effective and promising approach for tackling such general optimization problems is stochastic search. This refers to a collection of methods that use some sort of randomized mechanism to generate a sequence of iterates, e.g., candidate solutions, and then use the sequence of iterates to successively approximate the optimal solution. Over the past years, various stochastic search algorithms have been proposed in literature. These include approaches such as simulated annealing \cite{kirkpatrick:1983}, genetic algorithms \cite{goldberg:1989}, tabu search \cite{glover:1990}, pure adaptive search \cite{zabinsky:2003}, and sequential Monte Carlo simulated annealing \cite{zhou:2011}, which produce a sequence of candidate solutions that are gradually improving in performance; the nested partitions method \cite{shi:2000}, which uses a sequence of partitions of the feasible region as intermediate constructions to find high quality solutions; and the more recent class of model-based algorithms (see a survey by \cite{zlochin:2004a}), which construct a sequence of distribution models to characterize promising regions of the solution space.

This paper focuses on model-based algorithms. These algorithms typically assume a sampling distribution (i.e., a probabilistic model), often within a parameterized family of distributions, over the solution space, and iteratively carry out the two interrelated steps: (1) draw candidate solutions from the sampling distribution; (2) use the evaluations of these candidate solutions to update the sampling distribution. The hope is that  at every iteration the sampling distribution is biased towards the more promising regions of the solution space, and will eventually concentrate on one or more of the optimal solutions. Examples of model-based algorithms include ant colony optimization \cite{dorigo:1997, dorigo:2005}, annealing adaptive search (AAS) \cite{Romeijn:1994a}, probability collectives (PCs) \cite{wolpert:2004}, the estimation of distribution algorithms (EDAs) \cite{larranaga:1999a, muhlenbein:1996a}, the cross-entropy (CE) method \cite{rubinstein:2001a}, model reference adaptive search (MRAS) \cite{hu:2007a}, and the interacting-particle algorithm \cite{molvalioglu:2009, molvalioglu:2010}. The various model-based algorithms mainly differ in their ways of updating the sampling distribution. Recently, \cite{hu:2012a} showed that the updating schemes in some model-based algorithms can be viewed under a unified framework. The basic idea is to
convert the original optimization problem into a sequence of stochastic optimization problems with differentiable structures, so that the distribution updating schemes in these algorithms can be equivalently transformed into the form of stochastic approximation procedures for solving the sequence of stochastic optimization problems.

Because model-based algorithms work with a population of candidate solutions at each iteration, they demonstrate more robustness in exploring the solution space as compared with their classical counterparts that work with a single candidate solution each time (e.g., simulated annealing). The main motivation of this paper is to integrate this robustness feature of model-based algorithms into familiar gradient-based tools from classical differentiable optimization to facilitate the search for good sampling distributions.
The underlying idea is to reformulate the original (possibly non-differentiable) optimization problem into a {\em differentiable} optimization problem over the parameter space of the sampling distribution, and then use a direct gradient search method on the parameter space to solve the new formulation.
This leads to a natural algorithmic framework that combines the advantages of both methods: the fast convergence of gradient-based methods and the global exploration of stochastic search. Specifically, each iteration of our proposed method consists of the following two steps: (1) generate candidate solutions from the current sampling distribution; (2) update the parameters of the sampling distribution using a direct gradient search method. Although there are a variety of gradient-based algorithms that are applicable in step (2) above, in this paper we focus on a particular algorithm that uses a quasi-Newton like procedure to update the sampling distribution parameters. Note that since the algorithm uses only the information contained in the sampled solutions, it differs from the quasi-Newton method in deterministic optimization, in that there is an extra Monte Carlo sampling noise involved at each parameter updating step. We show that this stochastic version of quasi-Newton iteration can be expressed in the form of a generalized Robbins-Monro algorithm, and this in turn allows us to use the existing tools from stochastic approximation theory to analyze the asymptotic convergence and convergence rate of the proposed algorithm.

The rest of the paper is organized as follows. We introduce the problem setting formally in Section~\ref{sec:2}. Section~\ref{sec:3} provides a description of the proposed algorithm along with the detailed derivation steps. In Section~\ref{sec:4}, we analyze the asymptotic properties of the algorithm, including both convergence and convergence rate.
Some preliminary numerical study are carried out in Section~\ref{s:numerical} to illustrate the performance of the algorithm. Finally, we conclude this paper in Section~\ref{sec:6}. All the proofs are contained in the Appendix.

\section{Problem Formulation}\label{sec:2}
Consider the maximization problem
\begin{equation} \label{max H}
x^* \in \arg\max_{x \in \mathcal{X}} {H(x)}, ~~ \mathcal{X} \subseteq \mathbb{R}^n.
\end{equation}
where the solution space $\mathcal{X}$ is a nonempty compact set in $\mathbb{R}^n$, and $H: \mathcal{X}\rightarrow \mathbb{R}$ is a real-valued
function. Denote the optimal function value as $H^*$, i.e., there exists an $x^*$ such that $H(x)\leq H^* \triangleq H(x^*)$, $\forall
x\in\mathcal{X}$. Assume that $H$ is bounded on $\mathcal{X}$, i.e., $\exists H_{lb} > -\infty,~H_{ub} <
\infty$ s.t. $H_{lb} < H(x)< H_{ub}$, $\forall x\in \mathcal{X}$. We consider problems where the objective function $H(x)$ lacks nice structural properties such as differentiability and convexity and could have multiple local optima.

Motivated by the idea of using a sampling distribution/probabilistic model in model-based optimization, we let $\{f(x; \theta)| \theta \in \Theta \subseteq \mathbb{R}^{d} \}$ be a parameterized family of probability density functions (pdfs) on $\mathcal{X}$ with $\Theta$ being a parameter space. Intuitively, this collection of pdfs can be viewed abstractly as probability models characterizing our knowledge or belief of the promising regions of the solution space.
It is easy to see that
$$
\int{H(x)f(x; \theta)dx} \leqslant H^*,~~\forall \theta \in \Theta.
$$
In this paper, we simply write $\int$ with the understanding that the integrals are taken over $\mathcal{X}$. Note that the equality on the righthand side above is achieved whenever there exists an optimal parameter under which the parameterized probability distribution will assign all of its probability mass to a subset of the set of global optima of (\ref{max H}). Hence, one natural idea to
solving (\ref{max H}) is to transform the original problem
into an expectation of the objective function under the parameterized distribution and try to
 find the best parameter $\theta^*$ within the parameter space $\Theta$ such that the expectation under $f(x,\theta^*)$ can be made as large as possible, i.e.,
\begin{equation} \label{max Hf}
\theta^* = \arg\max_{\theta \in \Theta}{\int{H(x)f(x; \theta)dx}}.
\end{equation}
So instead of considering directly the original function $H(x)$ that is possibly non-differentiable and discontinuous in $x$, we now consider the new objective function $\int{H(x)f(x; \theta)dx}$ that is continuous on the parameter space and usually differentiable with respect to $\theta$. For example, under mild conditions the differentiation can be brought into the integration to apply on the p.d.f. $f(x;\theta)$, which is differentiable given an appropriate choice of the distribution family such as an exponential family of distributions. 

The formulation of (\ref{max Hf}) suggests a natural integration of stochastic search methods on the solution space $\mathcal{X}$ with gradient-based optimization techniques on the continuous parameter space. Conceptually, that is to iteratively carry out the following two steps:
\begin{enumerate}
  \item Generate candidate solutions from $f(x;\theta)$ on the solution space $\mathcal{X}$.
  \item Use a gradient-based method for the problem (\ref{max Hf}) to update the parameter $\theta$.
\end{enumerate}
The motivation is to speed up stochastic search with a guidance on the parameter space, and hence combine the advantages of both methods: the fast convergence of gradient-based methods and the global exploration of stochastic search methods.  Even though problem (\ref{max Hf}) may be non-concave and multi-modal in $\theta$, the sampling from the entire original space $\mathcal{X}$ compensates the local exploitation along the gradient on the parameter space. In fact, our algorithm developed later will automatically adjust the magnitude of the gradient step on the parameter space according to the global information, i.e., our belief about the promising regions of the solution space.

For algorithmic development later, we introduce a shape function $S_{\theta}: \mathbb{R} \rightarrow \mathbb{R}^+$, where the subscript $\theta$ signifies the possible dependence of the shape function on the parameter $\theta$. The function $S_{\theta}$ satisfies the following conditions:
\begin{enumerate}
\item[(a)] For every $\theta$, $S_{\theta}(y)$ is nondecreasing in $y$ and bounded from above and below for bounded $y$, with the lower bound being away from zero. Moreover, for every fixed $y$, $S_{\theta}(y)$ is continuous in $\theta$;
\item[(b)] The set of optimal solutions $\{\arg \max_{x \in \mathcal{X}} S_{\theta}(H(x))\}$ is a non-empty subset of $\{\arg \max_{x \in \mathcal{X}} H(x)\}$, the set of optimal solutions of the original problem (\ref{max H}).
\end{enumerate}
Therefore, solving (\ref{max H}) is equivalent to solving the following problem
\begin{equation} \label{max UH}
x^* \in \arg \max_{x \in \mathcal{X}} {S_{\theta}(H(x))}.
\end{equation}
The main reason of introducing the shape function $S_{\theta}$ is to ensure positivity of the objective function $S_{\theta}(H(x))$ under consideration, since  $S_{\theta}(H(x))$ will be used to form a probability density function later. Moreover, the choice of $S_{\theta}$ can also be used to adjust the trade-off between exploration and exploitation in stochastic search. One choice of such a shape function, similar to the level/indicator function used in the CE method and MRAS, is
\begin{equation} \label{quantile}
S_{\theta}(H(x))=(H(x)-H_{lb})\frac{1}{1+e^{-S_0(H(x)-\gamma_{\theta})}},
\end{equation}
where $S_0$ is a large positive constant, and $\gamma_{\theta}$ is the $(1-\rho)$-quantile
$$
\gamma_{\theta} \triangleq \sup_{l}\left\{l: P_{\theta}\{x\in \mathcal{X}: H(x) \geq l\} \geq \rho\right\},
$$
where $P_{\theta}$ denotes the probability with respect to $f(\cdot;\theta)$. Notice that $1/(1+e^{-S_0(H(x)-\gamma_{\theta})})$ is a continuous approximation of the indicator function $I\{H(x)\geq \gamma_{\theta}\}$, this shape function $S_{\theta}$ essentially prunes the level sets below $\gamma_{\theta}$. By varying $\rho$, we can adjust the percentile of elite samples that are selected to update the next sampling distribution: the bigger $\rho$, the less elite samples selected and hence more emphasis is put on exploiting the neighborhood of the current best solutions. Sometimes the function $S_{\theta}$ could also be chosen to be independent of $\theta$, i.e., $S_{\theta} = S: \mathbb{R} \rightarrow \mathbb{R}^+$, such as the function $S(y) = \exp(y)$.

For an arbitrary but fixed $\theta' \in \mathbb{R}^d$, define the function
$$
L(\theta; \theta') \triangleq \int{S_{\theta'}(H(x))f(x; \theta)dx}.
$$
According to the conditions on $S_{\theta}$, it always holds that
$$
0 < L(\theta; \theta') \leq S_{\theta'}(H^*)~~\forall\,\theta,
$$
and the equality is achieved if there exists an optimal parameter such that the probability mass of the parameterized distribution is concentrated only on a subset of the set of global optima. Following the same idea that leads to (\ref{max Hf}), solving (\ref{max UH}) and thus (\ref{max H}) can be converted to the problem of trying to find the best parameter $\theta^*$ within the parameter space by solving the following maximization problem:
\begin{equation} \label{max UHf}
\theta^* = \arg\max_{\theta \in \Theta}{L(\theta; \theta')}.
\end{equation}
Same as problem (\ref{max Hf}), $L(\theta; \theta')$ may be nonconcave and multi-modal in $\theta$.

\section{Gradient-Based Adaptive Stochastic Search}\label{sec:3}
Following the formulation in the previous section, we propose a stochastic search algorithm that carries out the following two steps at each iteration: let $\theta_k$ be the parameter obtained at the $k^{th}$ iteration,
\begin{enumerate}
  \item Generate candidate solutions from $f(x;\theta_k)$.
  \item Update the parameter to $\theta_{k+1}$ using a quasi Newton's iteration for $\max_{\theta}L(\theta;\theta_k)$.
\end{enumerate}
Assuming it is easy to draw samples from $f(x;\theta)$, then the main obstacle is to find expressions of the gradient and Hessian of $L(\theta;\theta_k)$ that can be nicely estimated using the samples from $f(x; \theta)$. To overcome this obstacle, we choose $\{f(x;\theta)\}$ to be an exponential family of densities defined as below.
\begin{Definition} \label{def:exp_fam}
A family $\{f(x;\theta): \theta \in \Theta\}$ is an exponential family of densities if it satisfies
\begin{equation} \label{exp family}
f(x;\theta) = \exp\{ \theta^T T(x) - \phi(\theta)\},  ~~\phi(\theta) = \ln \left\{\int \exp(\theta^T T(x)) dx \right\}.
\end{equation}
where $T(x) = [T_1(x), T_2(x), \ldots, T_{d}(x)]^T$ is the vector of sufficient statistics, $\theta  = [\theta_1, \theta_2, \ldots, \theta_{d}]^T$ is the vector of natural parameters, and $\Theta=\{\theta\in \mathbb{R}^d:\,|\phi(\theta)|<\infty \}$ is the natural parameter space with a nonempty interior.
\end{Definition}

Define the density function
\begin{equation} \label{p definition}
p(x; \theta) \triangleq \frac{S_{\theta}(H(x))f(x; \theta)}{\int S_{\theta}(H(x))f(x; \theta)dx}= \frac{S_{\theta}(H(x))f(x; \theta)}{L(\theta; \theta)}.
\end{equation}
With $f(\cdot; \theta)$ from an exponential family, we propose the following updating scheme for $\theta$ in  step 2 above:
\begin{equation}\label{theta ideal}
\theta_{k+1} = \theta_k + \alpha_k(\mathrm{Var}_{\theta_k}[T(X)] + \epsilon I)^{-1}\left(E_{p_{k}}[T(X)] - E_{\theta_k}[T(X)] \right),
\end{equation}
where $\epsilon > 0$ is a small positive number, $\alpha_k > 0$ is the step size, $E_{p_k}$ denotes the expectation with respect to $p(\cdot; \theta_k)$, and $E_{\theta_k}$ and $\mathrm{Var}_{\theta_k}$ denote the expectation and variance taken with respect to $f(\cdot; \theta_k)$, respectively.
The role of $\epsilon I$ is to ensure the positive definiteness of $(\mathrm{Var}_{\theta_k}[T(X)] + \epsilon I)$ such that it can be inverted. The term $(E_{p_k}[T(X)] - E_{\theta_k}[T(X)])$ is an ascent direction of $L(\theta; \theta_k)$, which will be shown in the next section.

To implement the updating scheme (\ref{theta ideal}), the term $E_{p_{k}}[T(X)]$ is often not analytically available and needs to be estimated. Suppose $\{x_1, \ldots, x_{N_k}\}$ are  independent and identically distributed (i.i.d.) samples drawn from  $f(x;\theta_k)$. Since
$$
E_{p_k}[T(X)] = E_{\theta_k}\left[T(X)\frac{p(X;\theta_k)}{f(X;\theta_k)}\right],
$$
we compute the weights $\{w_k^i\}$ for the samples $\{x^i_k\}$ according to
\begin{eqnarray*}
&& w^i_k \propto \frac{p(x^i_k; \theta_k)}{f(x^i_k; \theta_k)} \propto S_{\theta_k}(H(x^i_k)),~~i = 1, \ldots, N_k, \\
&& \sum_{i=1}^N{w_k^i}=1.
\end{eqnarray*}
Hence, $E_{p_k}[T(X)]$ can be approximated by
\begin{equation} \label{tildeEp}
\widetilde{E}_{p_k}[T(X)]  = \sum_{i=1}^{N_k}{w_k^iT(x_k^i)}.
\end{equation}
Some forms of the function $S_{\theta_k}(H(x))$ have to be approximated by samples as well. For example, if $S_{\theta_k}(H(x))$ takes the form (\ref{quantile}), the quantile $\gamma_{\theta_k}$  needs to be estimated by the sample quantile. In this case, we denote the approximation by $\widehat{S}_{\theta_k}(H(x))$, and evaluate the normalized weights according to
$$
\widehat{w}_i^k \propto \widehat{S}_{\theta_k}(H(x^i_k)), ~~i =1, \ldots, N_k.
$$
Then the term $E_{p_k}[T(X)]$ is approximated by
\begin{equation} \label{hatEp}
\widehat{E}_{p_k}[T(X)]  = \sum_{i=1}^{N_k}{\widehat{w}_k^iT(x_k^i)}.
\end{equation}
In practice, the variance term $\mathrm{Var}_{\theta_k}[T(X)]$ in (\ref{theta ideal}) may not be directly available or could be too complicated to compute analytically, so it also often needs to be estimated by samples:
\begin{eqnarray}
\widehat{\mathrm{Var}}_{\theta_k}[T(X)] &=& \frac{1}{N_k-1}\sum_{i=1}^{N_k}{T(x_k^i)T(x_k^i)^T} - \frac{1}{N_k^2-N_k} \left(\sum_{i=1}^{N_k}{T(x_k^i)}\right)\left(\sum_{i=1}^{N_k}{T(x_k^i)}\right)^T. \label{hessian theta appx}
\end{eqnarray}
The expectation term $E_{\theta_k}[T(X)]$ can be evaluated analytically in most cases. For example, if $\{f(\cdot;\theta_k)\}$ is chosen as the Gaussian family, then $E_{\theta_k}[T(X)]$ reduces to the mean and second moment of the Gaussian distribution.

Based on the updating scheme of $\theta$, we propose the following algorithm for solving (\ref{max H}).
\begin{algorithm}[H]\caption{\textbf{Gradient-Based Adaptive Stochastic Search (GASS)}} \label{alg: GASS2}
\begin{enumerate}
	\item {\em Initialization:} choose an exponential family of densities $\{f(\cdot; \theta)\}$, and specify a small positive constant $\epsilon$, initial parameter $\theta_0$, sample size sequence $\{N_k\}$, and step size sequence $\{\alpha_k\}$. Set $k=0$.
	\item {\em Sampling:} draw samples $x_{k}^i \iid f(x;\theta_k), i=1,2,\ldots,N_k$.
	\item {\em Estimation:} compute the normalized weights $\widehat{w}_k^i$ according to
	$$
	\widehat{w}_k^i = \frac{\widehat{S}_{\theta_k}(H(x_k^i))}{\sum_{j=1}^{N_k}{\widehat{S}_{\theta_k}(H(x_k^j))}},
	$$
	and then compute $\widehat{E}_{p_k}[T(X)]$ and $\widehat{\mathrm{Var}}_{\theta_k}[T(X)] $  respectively according to (\ref{hatEp}) and (\ref{hessian theta appx}).
\item {\em Updating:} update the parameter $\theta$ according to
    $$
    \ez{\theta_{k+1} = \Pi_{\tilde{\Theta}}\left\{ \theta_k + \alpha_k(\widehat{\mathrm{Var}}_{\theta_k}[T(X)] + \epsilon I)^{-1}(\widehat{E}_{p_k}[T(X)] - E_{\theta_k}[T(X)]) \right\},}
    $$
 where $\tilde{\Theta}\subseteq \Theta$ is a non-empty compact  connected constraint set, and $\Pi_{\tilde{\Theta}}$ denotes the projection operator that projects an iterate back onto the set $\tilde{\Theta}$ by choosing the closest point in $\tilde{\Theta}$.
\item {\em Stopping:} check if some stopping criterion is satisfied. If yes, stop and return the current best sampled solution; else, set $k := k+1$ and go back to step 2.
\end{enumerate}
\end{algorithm}

In the above algorithm, at the $k^{th}$  iteration candidate solutions are drawn from  the sampling distribution $f(\cdot; \theta_k)$, and then are used to estimate the quantities in the updating equation for $\theta_k$ so as to generate the next sampling distribution $f(\cdot; \theta_{k+1})$. To develop an intuitive understanding of the algorithm, we consider the special setting $T(X) = X$, in which case the term $\widehat{\mathrm{Var}}_{\theta_k}[T(X)]$ basically measures how widespread the candidate solutions are. Since the magnitude of the ascent step is determined by $(\widehat{\mathrm{Var}}_{\theta_k}[T(X)]+ \epsilon I)^{-1}$,  the algorithm takes smaller ascent steps to update $\theta$ when the candidate solutions are more widely spread (i.e., $\widehat{\mathrm{Var}}_{\theta_k}[X]$ is larger), and takes larger ascent steps when the candidate solutions are more concentrated (i.e.,  $\widehat{\mathrm{Var}}_{\theta_k}[X]$ is smaller). It means that exploitation of the local structure is adapted to our belief about the promising regions of the solution space: we will be more conservative in exploitation if we are uncertain about where the promising regions are and more progressive otherwise. Note that the projection operator at step $4$ is primarily used to ensure the numerical stability of the algorithm. It prevents the iterates of the algorithm from becoming too big in practice and ensures the sequence $\{\theta_k\}$ to stay bounded as the search proceeds. For simplicity, we will assume that $\tilde{\Theta}$ is a hyper-rectangle and takes the form
$\tilde\Theta=\{\theta\in\Theta:a_i\leq \theta_i \leq b_i \}$ for constants $a_i<b_i$, $i=1,\ldots,d$; other more general choices of $\tilde{\Theta}$ may also be used (see, e.g., Section 4.3 of \cite{kushner:2003}). Intuitively, such a constraint set should be chosen sufficiently large in practice so that the limits of the recursion at step 4 without the projection are contained in its interior.

\subsection{Accelerated GASS}

GASS can be viewed as a stochastic approximation (SA) algorithm, which we will show in more details in the next section. To improve the convergence rate of SA algorithms, \cite{polyak:1990} and \cite{ruppert:1991} first proposed to take the average of the $\theta$ values generated by previous iterations, which is often referred to as Polyak (or Polyak-Ruppert) averaging. The original Polyak averaging technique is ``offline'', i.e., the averages are not fed back into the iterates of $\theta$, and hence the averages are not useful for guiding the stochastic search in our context. However, there is a variation, Polyak averaging with online feedback (c.f. pp. 75 - 76 in \cite{kushner:2003}), which is not optimal as the original Polyak averaging but also enhances the convergence rate of SA. Using the Polyak averaging with online feedback, the parameter $\theta$ will be updated according to
\begin{equation} \label{theta avg2}
\ez{\theta_{k+1}=\Pi_{\tilde{\Theta}}\left\{\theta_k + \alpha_k\left(\widehat{\mathrm{Var}}_{\theta_k}[T(X)] + \epsilon I \right)^{-1}(\widehat{E}_{p_k}[T(X)] - E_{\theta_k}[T  (X)] ) + \alpha_kc(\bar{\theta}_k-\theta_k)\right\},}
\end{equation}
where the constant $c$ is the feedback weight, and $\bar{\theta}_k$ is the average
$$
\bar{\theta}_k=\frac{1}{k}\sum_{i=1}^k\theta_i,
$$
which  can be calculated recursively by
\begin{equation} \label{thetabar}
\bar{\theta}_{k}=\frac{k-1}{k}\bar{\theta}_{k-1}+\frac{\theta_{k}}{k}.
\end{equation}
With this parameter updating scheme, we propose the following algorithm.
\begin{algorithm}[H]\caption{\textbf{Gradient-based Adaptive Stochastic Search with Averaging (GASS\_avg)}}  \label{alg:GSSA}
Same as Algorithm~\ref{alg: GASS2} except in step 4 the parameter updating follows (\ref{theta avg2}) and (\ref{thetabar}).
\end{algorithm}

\subsection{Derivation}
In this subsection,  we explain the rationale behind the updating scheme (\ref{theta ideal}). We first derive the expressions of the gradient and Hessian of $L(\theta;\theta')$ as given below.
\begin{Proposition} \label{prop: L gradient and hessian}
Assume that $f(x;\theta)$ is twice differentiable on $\Theta$ and that  $\nabla_{\theta}f(x;\theta)$ and $\nabla^2_{\theta}f(x;\theta)$ are both bounded on $\mathcal{X}$ for any $\theta\in\Theta$. Then
\begin{eqnarray*}
\nabla_{\theta}L(\theta; \theta') &=& E_{\theta}[S_{\theta'}(H(X))\nabla_{\theta}\ln{f(X;\theta)}] \\
\nabla^2_{\theta}L(\theta; \theta') &=& E_{\theta}[S_{\theta'}(H(X))\nabla^2_{\theta}\ln{f(X;\theta)}] \nonumber \\
&& +~ E_{\theta}[S_{\theta'}(H(X))\nabla_{\theta}\ln{f(X;\theta)}\nabla_{\theta}\ln{f(X;\theta)}^T].
\end{eqnarray*}
Furthermore, if $f(x;\theta)$ is in an exponential family of densities defined by (\ref{exp family}), then the above expressions reduce to
\begin{eqnarray*}
\nabla_{\theta}{L(\theta; \theta')} &=& E_{\theta}[S_{\theta'}(H(X))T(X)] - E_{\theta}[S_{\theta'}(H(X))]E_{\theta}[T(X)], \label{L gradient}\\
\nabla_{\theta}^2{L(\theta; \theta')} &=& E_{\theta}\left[S_{\theta'}(H(X))(T(X)-E_{\theta}[T(X)])(T(X)-E_{\theta}[T(X)])^T\right]  \nonumber \\
 && -~ \mathrm{Var}_{\theta}[T(X)]E_{\theta}[S_{\theta'}(H(X))].
\end{eqnarray*}
\end{Proposition}


Notice that if we were to use Newton's method to update the parameter $\theta$, the Hessian $\nabla^2_{\theta}L(\theta; \theta')$ is not necessarily negative semidefinite to ensure the parameter updating is along the ascent direction of $L(\theta; \theta')$, so we need some stabilization scheme. One way is to approximate the Hessian by the second term on the righthand side with a small perturbation, i.e., $-(\mathrm{Var}_{\theta}[T(X)]+ \epsilon I)E_{\theta}[S_{\theta'}(H(X))]$, which is always negative definite. Thus, the parameter $\theta$ could be updated according to the following iteration
\begin{eqnarray}
\theta_{k+1} &=& \theta_k + \alpha_k \left((\mathrm{Var}_{\theta_k}[T(X)]+\epsilon I)E_{\theta_k}[S_{\theta_k}(H(X))]\right)^{-1} \nabla_{\theta}{L(\theta_k; \theta_k)}, \label{Ltheta}\\
&=& \theta_k + \alpha_k \left( \mathrm{Var}_{\theta_k}[T(X)] + \epsilon I \right)^{-1} \left( \frac{E_{\theta_k}[S_{\theta_k}(H(X))T(X)]}{E_{\theta_k}[S_{\theta_k}(H(X))]} - E_{\theta_k}[T(X)] \right), \nonumber
\end{eqnarray}
which immediately leads to the updating scheme (\ref{theta ideal}) given before. 

In the updating equation (\ref{theta ideal}), the term $E_{\theta_k}[S_{\theta_k}(H(X))]^{-1}$ is absorbed into $\nabla_{\theta}L(\theta_k;\theta_k)$, so we obtain a scale-free term $\left( E_{p_k}[T(X)] - E_{\theta_k}[T(X)] \right)$ that is not subject to the scaling of the function value of $S_{\theta_k}(H(x))$. It would be nice to have such a scale-free gradient so that we can employ other gradient-based methods more easily besides the above specific choice of a quasi-Newton method. Towards this direction, we consider a further transformation of the maximization problem (\ref{max UHf}) by letting
$$
l(\theta;\theta') = \ln{L(\theta;\theta')}.
$$
Since $\ln: R^+ \rightarrow R$ is a strictly increasing function, the maximization problem (\ref{max UHf}) is equivalent to
\begin{equation} \label{max l}
\theta^* = \arg \max_{\theta \in \mathbb{R}^d} {l(\theta;\theta')}.
\end{equation}
The gradient and the Hessian of  $l(\theta;\theta')$ are given in the following proposition.
\begin{Proposition} \label{prop: gradient and hessian}
Assume that $f(x;\theta)$ is twice differentiable on $\Theta$ and that  $\nabla_{\theta}f(x;\theta)$ and $\nabla^2_{\theta}f(x;\theta)$ are both bounded on $\mathcal{X}$ for any $\theta\in\Theta$. Then
\begin{eqnarray*}
\nabla_{\theta}{l(\theta; \theta')}|_{\theta = \theta'} &=& E_{p(\cdot; \theta')} [\nabla_{\theta}\ln{f(X;\theta')}]  \\
\nabla_{\theta}^2{l(\theta; \theta')}|_{\theta = \theta'} &=& E_{p(\cdot; \theta')}[\nabla_{\theta}^2 \ln{f(X;\theta')}]  + \mathrm{Var}_{p(\cdot ; \theta')}\left[ \nabla_{\theta}\ln{f(X;\theta')} \right].
\end{eqnarray*}
Furthermore, if $f(x;\theta)$ is in an exponential family of densities, then the above expressions reduce to
\begin{eqnarray*}
\nabla_{\theta}{l(\theta; \theta')}|_{\theta = \theta'} &=& E_{p(\cdot; \theta')}[T(X)] - E_{\theta'}[T(X)], \label{gradient theta}\\
\nabla_{\theta}^2{l(\theta; \theta')}|_{\theta = \theta'} &=& \mathrm{Var}_{p(\cdot; \theta')}[T(X)] - \mathrm{Var}_{\theta'}[T(X)].  \label{hessian theta}
\end{eqnarray*}
\end{Proposition}


Similarly as before, noticing that the Hessian $\nabla^2_{\theta}l(\theta';\theta')$ is not necessarily negative definite to ensure the parameter updating is along the ascent direction of $l(\theta;\theta')$, we approximate the Hessian by the slightly perturbed second term in $\nabla_{\theta}^2{l(\theta';\theta')}$, i.e., $-(\mathrm{Var}_{\theta'}[T(X)]+ \epsilon I)$. Then by setting
$$
\theta_{k+1} = \theta_k + \alpha_k \left(\mathrm{Var}_{\theta_k}[T(X)] + \epsilon I \right)^{-1}\nabla_{\theta}{l(\theta_k;\theta_k)},
$$
we again obtain exactly the same updating equation (\ref{theta ideal}) for $\theta$. The difference from (\ref{L gradient}) is that the gradient $\nabla_{\theta}l(\theta;\theta')$ is a scale-free term, and hence can be used in other gradient-based methods with easier choices of the step size. From the algorithmic viewpoint, it is better to consider the optimization problem (\ref{max l}) on $l(\theta;\theta')$ instead of the problem (\ref{max UHf}) on $L(\theta;\theta')$, even though both have the same global optima.

Although there are many ways to determine the positive definite matrix in front of the gradient in a quasi-Newton method, our choice of $\left(\mathrm{Var}_{\theta_k}[T(X)] + \epsilon I \right)^{-1}$ is not arbitrary but based on some principle. Without considering the numerical stability and thus dropping the term $\epsilon I$, the term $\mathrm{Var}_{\theta}[T(X)] = E[\nabla_{\theta}\ln{f(X;\theta)}(\nabla_{\theta}\ln{f(X;\theta)})^T] = E[-\nabla_{\theta}^2\ln{f(X;\theta)}]$ is the Fisher information matrix, whose inverse provides a lower bound on the variance of an unbiased estimator of the parameter $\theta$ (\cite{rao:1945}), leading to the fact that $(\mathrm{Var}_{\theta}[T(X)])^{-1}$ is the minimum-variance step size in stochastic approximation (\cite{meyn:2009}). Moreover, from the optimization perspective, the term $(\mathrm{Var}_{\theta}[T(X)])^{-1}$ relates the gradient search on the parameter space with the stochastic search on the solution space, and thus  adaptively adjusts the updating of the sampling distribution to our belief about the promising regions of the solution space. To see this more easily, let us consider $T(X) = X$.
Then $(\mathrm{Var}_{\theta}[X])^{-1}$ is smaller (i.e., the gradient step in updating $\theta$ is smaller) when the current sampling distribution is more flat, signifying the exploration of the solution space is still active and we do not have a strong belief (i.e. $f(\cdot;\theta)$) about promising regions; $(\mathrm{Var}_{\theta}[X])^{-1}$ is larger (i.e., the gradient step in updating $\theta$ is larger) when our belief $f(\cdot;\theta)$ is more focused on some promising regions.


\section{Convergence Analysis}\label{sec:4}
We will analyze the convergence properties of  GASS, resorting to methods and results in stochastic approximation (e.g., \cite{kushner:1978, kushner:2003, borkar:2008}). In GASS, $\nabla_{\theta}l(\theta;\theta_k)|_{\theta = \theta_k}$ is estimated by
\begin{eqnarray}
\widehat{\nabla}_\theta l(\theta_k; \theta_k) = \widehat{E}_{p_k}[T(X)] - E_{\theta_k}[T(X)]. \label{gradient theta appx1}
\end{eqnarray}
To simplify notations, we denote
\begin{eqnarray*}
\widehat{V}_k &\triangleq& \widehat{\mathrm{Var}}_{\theta_k}[T(X)]+\epsilon I,  ~~~ V_k \triangleq \mathrm{Var}_{\theta_k}[T(X)]+\epsilon I.
\end{eqnarray*}
Hence, the parameter updating iteration in GASS is
\begin{equation}\label{theta generation}
\theta_{k+1} = \Pi_{\tilde{\Theta}}\left\{\theta_k + \alpha_k \widehat{V}_k^{-1} \widehat{\nabla}_{\theta}l(\theta_k; \theta_k)\right\},
\end{equation}
which can be  rewritten in the form of a generalized Robbins-Monro algorithm
\begin{eqnarray}
\theta_{k+1} = \theta_k + \alpha_k[D(\theta_k) + b_k + \xi_k + z_k], \label{robbins-monro}
\end{eqnarray}
where
\begin{eqnarray*}
D(\theta_k) &=& \left( \mathrm{Var}_{\theta_k}[T(X)]+\epsilon I\right)^{-1}\nabla_{\theta}l(\theta_k; \theta_k), \\
b_k &=&  \widehat{V}_k^{-1}\left(\widehat{E}_{p_k}[T(X)] - \widetilde{E}_{p_k}[T(X)]\right), \nonumber \\
\xi_k &=& \left(\widehat{V}_k^{-1}- V_k^{-1}\right)\left(\widetilde{E}_{p_k}[T(X)]-E_{\theta_k}[T(X)]\right) + V_k^{-1}\left(\widetilde{E}_{p_k}[T(X)] - E_{p_k}[T(X)] \right), \nonumber
\end{eqnarray*}
and $z_k$ is the projection term satisfying
$\alpha_k z_k=\theta_{k+1}-\theta_k-\alpha_k[D(\theta_k)+b_k+\xi_k]$, the minimum Euclidean length vector that takes the current iterate back onto the constraint set. The term $D(\theta_k)$ is the gradient vector field,  $b_k$ is the bias due to the inexact evaluation of the shape function in $\widehat{E}_{p_k}[T(X)]$ ($b_k$ is zero if the shape function can be evaluated exactly), and $\xi_k$ is the noise term due to Monte Carlo sampling in the approximations $\widehat{\mathrm{Var}}_{\theta_k}[T(X)]$ and $\widetilde{E}_{p_k}[T(X)]$.

For a given $\theta\in\tilde\Theta$, we define a set $C(\theta)$ as follows: if $\theta$ lies in the interior of $\tilde\Theta$, let $C(\theta)=\{0\}$; if $\theta$ lies on the boundary of $\tilde\Theta$, define $C(\theta)$ as the infinite convex cone generated by the outer normals at $\theta$ of the faces on which $\theta$ lies (\cite{kushner:2003} pp. 106).
The difference equation (\ref{robbins-monro}) can be viewed as a noisy discretization of the constrained ordinary differential equation (ODE)
\begin{equation} \label{theta ode}
\dot{\theta}_t = D(\theta_t) + z_t,~~ z_t \in -C(\theta_t),~~~t \geq 0,
\end{equation}
where $z_t$ is the minimum force needed to keep the trajectory of the ODE in $\tilde{\Theta}$. Thus, the sequence of $\{\theta_k\}$ generated by (\ref{robbins-monro}) can be shown to asymptotically approach the solution set of the above ODE (\ref{theta ode}) by using the well-known ODE method. Let $\|\cdot\|$ denote the vector supremum norm (i.e., $\|x\| = \max\{|x_i|\}$) or the matrix max norm (i.e., $\|A\| = \max\{|a_{ij}|\}$). Let $\|\cdot\|_2$ denote the vector 2-norm (i.e., $\|x\|=\sqrt{x_1^2 + \ldots + x_n^2}$) or the matrix norm induced by the vector 2-norm (also called spectral norm for a square matrix, i.e., $\|A\|_2 = \sqrt{\lambda_{max}(A^*A)}$, where $A^*$ is the conjugate transpose of $A$ and $\lambda_{max}$ returns the largest eigenvalue).

To proceed to the formal analysis, we introduce the following notations and assumptions. We denote the sequence of increasing sigma-fields generated by all the samples up to the $k^{th}$ iteration  by
$$
\left\{\mathcal{F}_k = \sigma\left(\{x_0^i\}_{i=1}^{N_0}, \{x_1^i\}_{i=1}^{N_1}, \ldots, \{x_k^i\}_{i=1}^{N_k}\right), k = 0, 1, \ldots\right\}.
$$
Define notations
\begin{align*}
\bar {\mathbb{U}}_k:=\frac{1}{N_k}\sum_{i=1}^{N_k}\widehat S_{\theta_k}(H(x^i_k))T(x^i_k),~ & \bar {\mathbb{V}}_k:=\frac{1}{N_k}\sum_{i=1}^{N_k}\widehat S_{\theta_k}(H(x^i_k))\\
\tilde {\mathbb{U}}_k:=\frac{1}{N_k}\sum_{i=1}^{N_k}S_{\theta_k}(H(x^i_k))T(x^i_k),~ & \tilde {\mathbb{V}}_k:=\frac{1}{N_k}\sum_{i=1}^{N_k}S_{\theta_k}(H(x^i_k))\\
{\mathbb{U}}_k:=E_{\theta_k}[S_{\theta_k}(H(X))T(X)],~ & \mathbb{V}_k:=E_{\theta_k}[S_{\theta_k}(H(X))].
\end{align*}
\begin{Assumption} \ \\ \label{assump: gain samplesize}
(i) The step size sequence $\{\alpha_k\}$ satisfies $\alpha_k > 0$ for all $k$, $\alpha_k \searrow 0$ as $k \rightarrow \infty$, and $\sum_{k=0}^{\infty}{\alpha_k} = \infty$. \\
(ii) The sample size $N_k=N_0k^{\zeta}$, where $\zeta>0$; moreover, $\{\alpha_k\}$ and $\{N_k\}$ jointly satisfies $\frac{\alpha_k}{\sqrt{N_k}} =  O(k^{-\beta})$ for some constant $\beta > 1$. \\
(iii) The function $x \mapsto T(x)$ is bounded on $\mathcal{X}$. \\
(iv) For any $x$, $|\widehat{S}_{\theta_k}(H(x))- S_{\theta_k}(H(x))|\rightarrow 0$ w.p.1 as $N_k \rightarrow \infty$.
\end{Assumption}
In the above assumption, (i) is a typical assumption on the step size sequence in SA, which means that $\alpha_k$ diminishes not too fast. Assumption~\ref{assump: gain samplesize}(ii) provides a guideline on how to choose the sample size given a choice of the step size sequence, and shows that the sample size has to increase to infinity no slower than a certain speed. For example, if we choose $\alpha_k = \alpha_0k^{-\alpha}$ with $0 <\alpha <1$, then it is sufficient to choose $N_k = O(k^{2(\beta -\alpha)})$.
Assumption~\ref{assump: gain samplesize}(iii) holds true for many exponential families used in practice.
Assumption~\ref{assump: gain samplesize}(iv) is a sufficient condition to ensure the strong consistency of estimates, and is satisfied by many choices of the shape function $S_{\theta}$. For example, it is trivially satisfied if $S_{\theta} = S$, since $S(H(x))$ can be evaluated exactly for each $x$. If $S_{\theta}$ takes the form of (\ref{quantile}), Assumption~\ref{assump: gain samplesize}(iv) is also satisfied, as shown in the following lemma.

\begin{Lemma} \label{lem:quant}
Suppose the shape function takes the form
\begin{equation*}
S_{\theta_k}(H(x))=(H(x)-H_{lb})\frac{1}{1+e^{S_0(H(x)- \gamma_{\theta_k})}},
\end{equation*}
where $\gamma_{\theta_k} \triangleq \sup_{l}\left\{l: P_{\theta_k}\{x\in \mathcal{X}: H(x) \geq l\} \geq \rho\right\}$ is the unique $(1-\rho)$-quantile with respect to $f(\cdot; \theta_k)$. If $S_{\theta_k}(H(x))$ is estimated by $\widehat{S}_{\theta_k}(H(x))$
with the true quantile $\gamma_{\theta_k}$ being replaced by the sample $(1-\rho)$-quantile
$
\widehat{\gamma}_{\theta_k}  = H_{(\lceil (1-\rho)N_k\rceil)},
$
where $\lceil a \rceil$ is the smallest integer greater than $a$, and $H_{(i)}$ is the $i^{th}$ order statistic of the sequence $\{H(x_k^i), i=1, \ldots, N_k\}$. Then under the condition $N_k=\Theta(k^{\zeta})$ $\zeta>0$, we have that for every $x$,
$
\big|\widehat{S}_{\theta_k}(H(x))- S_{\theta_k}(H(x))\big|\rightarrow 0~w.p.1 ~~{\rm as}~ k \rightarrow \infty.
$
\end{Lemma}


The next lemma shows that the summed tail error goes to zero w.p.1.
\begin{Lemma}\label{lem:1}
Under Assumption~\ref{assump: gain samplesize}~(i)-(iii), for any $T>0$,
$$
\lim_{k\rightarrow \infty}\left\{\sup_{\{n:0\leq \sum_{i=k}^{n-1}\alpha_i\leq T\}}\left\|\sum_{i=k}^{n}{\alpha_i\xi_i}\right\| \right\} = 0,~~w.p.1.
$$
\end{Lemma}

Theorem~\ref{thm:SA} below shows that GASS generates a sequence $\{\theta_k\}$ that asymptotically approaches the limiting solution of the ODE (\ref{theta ode}) under the regularity conditions specified in Assumption~\ref{assump: gain samplesize}.
\begin{Theorem} \label{thm:SA}
Assume that $D(\theta_t)$ is continuous with a unique integral curve (i.e., the ODE (\ref{theta ode}) has a unique solution $\theta(t)$) and Assumption~\ref{assump: gain samplesize} holds.  Then the sequence $\{\theta_k\}$ generated by (\ref{theta generation}) converges to a \ez{limit set} of (\ref{theta ode}) w.p.1. Furthermore, if the \ez{limit} sets of (\ref{theta ode}) are isolated equilibrium points, then w.p.1 $\{\theta_k\}$ converges to a unique equilibrium point.
\end{Theorem}

For a given distribution family, Theorem~1 shows that our algorithm will identify a local/global optimal sampling distribution within the given family that
provides the best capability in generating an optimal solution to (\ref{max H}). From the viewpoint of maximizing $E_{\theta}[H(X)]$, the average function value under our belief of where promising solutions are located (i.e., the parameterized distribution $f(x,\theta)$), the convergence of the algorithm to a local/global optimum in the parameter space essentially gives us a local/global optimum of our belief about the function value.

\subsection{Asymptotic Normality of GASS}
In this section, we study the asymptotic convergence rate of Algorithm $1$ under the assumption that the parameter sequence $\{\theta_k\}$ converges to a unique equilibrium point $\theta^*$ of the ODE (\ref{theta ode}) in the interior of $\tilde\Theta$. This indicates that
there exists a small open neighborhood $\mathcal{N}(\theta^*)$ of $\theta^*$ such that the sequence $\{\theta_k\}$ will be contained in $\mathcal{N}(\theta^*)$ for $k$ sufficiently large w.p.1. Thus, the projection operator in (\ref{theta generation}) and $z_k$ in (\ref{robbins-monro})
can be dropped in the analysis, because the projected recursion will behave identically to an unconstrained algorithm in the long run.
Define $\mathcal{L}(\theta)=\nabla_{\theta'}l(\theta';\theta)|_{\theta'=\theta}$ and let $J_{\mathcal{L}}$ be the Jacobian of $\mathcal{L}$. Under our conditions, it immediately follows from (\ref{theta ode}) that $C(\theta^*)=\{0\}$ and $\mathcal{L}(\theta^*)=0$.
Since $\mathcal{L}$ is the gradient of some underlying function $F(\theta)$, $J_{\mathcal{L}}$ is the Hessian of $F$ and Algorithm 1 is essentially a gradient-based algorithm for maximizing $F(\theta)$. Therefore, it is reasonable to expect that the following assumption holds:
\begin{Assumption}
The Hessian matrix $J_{\mathcal{L}}(\theta)$ is continuous and symmetric negative definite in the neighborhood $\mathcal{N}(\theta^*)$ of $\theta^*$.
\end{Assumption}
We consider a standard gain sequence $\alpha_k=\alpha_0/k^{\alpha}$ for constants $\alpha_0>0$ and $0<\alpha<1$, a polynomially increasing sample size \ez{$N_k=N_0 k^{\zeta}$} with $N_0\geq 1$ and $\zeta>0$.

By dropping the projection operator in (\ref{theta generation}), we can rewrite the equation in the form:
$$\delta_{k+1}=\delta_k+k^{-\alpha}\Phi_k\mathcal{L}(\theta_k)+k^{-\alpha}\Phi_k\Big( \frac{\bar {\mathbb{U}}_k}{\bar {\mathbb{V}}_k}-\frac{{\mathbb{U}}_k}{{\mathbb{V}}_k}\Big),$$
where $\delta_k=\theta_k-\theta^*$ and $\Phi_k=\alpha_0(\widehat{\mathrm{Var}}_{\theta_k}(T(X))+\epsilon I)^{-1}$. Next, by using a first order Taylor expansion of $\mathcal{L}(\theta_k)$ around the neighborhood of $\theta^*$ and the fact that $\mathcal{L}(\theta^*)=0$, we have
$$\delta_{k+1}=\delta_k+ k^{-\alpha}\Phi_k J_{\mathcal{L}}(\tilde\theta_k)\delta_k+
 k^{-\alpha} \Phi_k \Big(\frac{\bar {\mathbb{U}}_k}{\bar {\mathbb{V}}_k}-\frac{{\mathbb{U}}_k}{{\mathbb{V}}_k} \Big),$$
where $\tilde \theta_k$ lies on the line segment from $\theta_k$ to $\theta^*$. For a given positive constant $\tau>0$, the above equation can be further written in the form of a recursion in \cite{fabian:1968}:
$$\delta_{k+1}=(I-k^{-\alpha}\Gamma_k) \delta_k+k^{-(\alpha+\tau)/2} \Phi_kW_k+k^{-\alpha-\tau/2}T_k,$$
where $\Gamma_k=-\Phi_k J_{\mathcal{L}}(\tilde \theta_k)$, $W_k= k^{(\tau-\alpha)/2}\big(\frac{\tilde {\mathbb{U}}_k}{\tilde {\mathbb{V}}_k}-E_{\theta_k}\big[\frac{\tilde {\mathbb{U}}_k}{\tilde {\mathbb{V}}_k}\big|\mathcal{F}_{k-1}\big] \big)$, and $T_k= k^{\tau/2}\Phi_k\big(\frac{\bar {\mathbb{U}}_k}{\bar {\mathbb{V}}_k}-\frac{\tilde {\mathbb{U}}_k}{\tilde {\mathbb{V}}_k}+E_{\theta_k}\big[\frac{\tilde {\mathbb{U}}_k}{\tilde {\mathbb{V}}_k}\big|\mathcal{F}_{k-1}\big]-\frac{{\mathbb{U}}_k}{{\mathbb{V}}_k}  \big)$.
The basic idea of the rate analysis is to show that the sequence of amplified differences $\{k^{\tau/2}\delta_k\}$ converges in distribution to a normal random variable with mean zero and constant covariance matrix. To this end, we show that all sufficient conditions in Theorem 2.2 in \cite{fabian:1968} are satisfied in our setting. We begin with a strengthened version of Assumption 1($iv$).

\begin{Assumption}\ \\
For a given constant $\tau>0$ and $x\in\mathcal{X}$, $k^{\tau/2}|\widehat S_{\theta_k}(H(x))-S_{\theta_k}(H(x))|\rightarrow 0$ as $k\rightarrow \infty$ w.p.1.
\end{Assumption}

Assumption 3 holds trivially when $S_{\theta}$ is a deterministic function that is independent of $\theta$. In addition, if sample quantiles are involved in the shape function and $S_{\theta_k}(H(x))$ takes the form (\ref{quantile}), then the assumption can also be justified under some additional mild regularity conditions; cf. e.g., \cite{hu:2012a}.

Let $\Phi=\alpha_0({\mathrm{Var}}_{\theta^*}(T(X))+\epsilon I)^{-1}$ and $\Gamma=-\Phi J_{\mathcal{L}}(\theta^*)$. The following result shows condition (2.2.1) in Theorem 2.2 of \cite{fabian:1968}.
\begin{Lemma}\label{lem1}
Assume Assumptions 1 and 2 hold, we have $\Phi_k\rightarrow \Phi$ and $\Gamma_k \rightarrow \Gamma$ as $k \rightarrow \infty$ w.p.1. In addition, if Assumption 1(iv) is replaced with Assumption 3 and \ez{$N_k=N_0k^{\zeta}$} with $\zeta>\tau/2$, then $T_k \rightarrow 0$ as $k \rightarrow \infty$ w.p.1.
\end{Lemma}

In addition, the noise term $W_k$ has the following property, which justifies condition 2.2.2 in \cite{fabian:1968}.
\begin{Lemma}\label{lem2}
$E_{\theta_k}[W_k|\mathcal{F}_{k-1}]=0$. In addition, let $\tau$ be a given constant satisfying $\tau> \alpha$. If Assumption 1 holds and \ez{$N_k=N_0k^{\tau-\alpha}$}, then there exists a positive semi-definite matrix $\Sigma$ such that $\lim_{k\rightarrow \infty}E_{\theta_k}[W_k W_k^T|\mathcal{F}_{k-1}]=\Sigma$ w.p.1, and $\lim_{k\rightarrow \infty}E[I\{\|W_k\|^2\geq r k^{\alpha}\}\|W_k \|^2]=0$ $\forall r>0$.
\end{Lemma}

The following asymptotic normality results then follows directly from Theorem 2.2 in \cite{fabian:1968}.
\begin{Theorem}
Let $\alpha_k=\alpha_0/k^{\alpha}$ for $0<\alpha<1$. For a given constant $\tau>2\alpha$, let $N_k=N_0k^{\tau-\alpha})$. Assume the convergence of the sequence $\{\theta_k\}$ occurs to a unique equilibrium point $\theta^*$ w.p.1. If Assumptions 1, 2, and 3 hold, then
$$k^{\frac{\tau}{2}}(\theta_k-\theta^*)\xrightarrow{~dist~} N(0,Q\mathcal{M}Q^T),$$
where $Q$ is an orthogonal matrix such that $Q^T(-J_{\mathcal{L}}(\theta^*))Q=\Lambda$ with $\Lambda$ being a diagonal matrix, and the $(i,j)^{th}$ entry of the matrix $\mathcal{M}$ is given by
$\mathcal{M}_{(i,j)}=(Q^T\Phi\Sigma \Phi^T Q)_{(i,j)}(\Lambda_{(i,i)}+\Lambda_{(j,j)})^{-1}$.
\end{Theorem}
Theorem 2 shows the asymptotic rate at which the noise caused by Monte-Carlo random sampling in GASS will be damped out as the number of iterations $k\rightarrow \infty$. This rate, as indicated in the theorem, is on the order of $O(1/\sqrt{k^{\tau}})$. This implies that the noise can be damped out arbitrarily fast by using a sample size sequence $\{N_k\}$ that increases sufficiently fast as $k\rightarrow \infty$. However, we note that this rate result is stated in terms of the number of iterations $k$, not the sample size $N_k$. Therefore, in  practice, there is the need to carefully balance the tradeoff between the choice of large values of $N_k$ to increase the algorithms's asymptotic rate and the use of small values of $N_k$ to reduce the per iteration computational cost.

\section{Numerical Experiments} \label{s:numerical}

We test the proposed algorithms GASS, GASS\_avg on some benchmark continuous optimization problems selected from
\cite{hu:2007a} and \cite{hu:2012a}. To fit in the maximization framework where our algorithms are proposed, we take the
negative of those objective functions that are originally for minimization problems. The ten benchmark problems are listed as
below.

\begin{enumerate}
\item[(1)] Dejong's 5th function (n=2, $-50\leq x_i\leq 50$)
\begin{equation*}\label{Dejong5}
H_1(x)=-\left[0.002+\sum_{j=1}^{25}\frac{1}{j+\sum_{i=1}^2(x_i-a_{ji})^6}\right]^{-1},
\end{equation*}
where $a_{j1}=(-32, -16, 0, 16, 32, -32, -16, 0, 16, 32, -32, -16,
0, 16, 32, -32, -16, 0, 16,\\
 32, -32, -16, 0, 16, 32)$ and
$a_{j2}=(-32, -32, -32, -32, -32, -16, -16, -16, -16, -16, 0,\\
 0, 0, 0, 0, 16, 16, 16, 16, 16, 32, 32, 32, 32, 32)$. The global
optimum is at $x^*=(-32,-32)^T$, and $H^*\approx-0.998$.
\item[(2)] Shekel's function (n=4, $0\leq x_i\leq 10$ )
\begin{equation*}\label{Shekel}
H_2(x)=\sum_{i=1}^{5}\left((x-a_i)^T(x-a_i)+c_i\right)^{-1},
\end{equation*}
where $a_1=(4,4,4,4)^T$, $a_2=(1, 1, 1, 1)^T$, $a_3=(8, 8, 8, 8)^T$,
$a_4=(6, 6, 6, 6)^T$, $a_5=(3, 7, 3, 7)^T$, and $c=(0.1, 0.2, 0.2,
0.4, 0.4)$. $x^*=(4, 4, 4, 4)^T$, $H^*\approx10.153$.
\item[(3)] Powel singular function (n=50, $-50\leq x_i\leq 50$)
\begin{equation*}\label{Powel}
H_3(x)=-\sum_{i=2}^{n-2}\left[(x_{i-1}+10x_i)^2+5(x_{i+1}-x_{i+2})^2+(x_{i}-2x_{i+1})^4+10(x_{i-1}-x_{i+2})^4\right]-1,
\end{equation*}
where $x^*=(0,\cdots,0)^T$, $H^*=-1$.
\item[(4)] Rosenbrock function (n=10, $-10\leq x_i\leq 10$)
\begin{equation*}\label{Rosenbrock}
H_4(x)=-\sum_{i=1}^{n-1}\left[100(x_{i+1}-x_i^2)^2+(x_i-1)^2\right]-1,
\end{equation*}
where $x^*=(1,\cdots,1)^T$, $H^*=-1$.
\item[(5)] Griewank function (n=50, $-50\leq x_i\leq 50$)
\begin{equation*}\label{Griewank}
H_5(x)=-\frac{1}{4000}\sum_{i=1}^{n}x_i^2+\prod_{i=1}^{n}\cos\left(\frac{x_i}{\sqrt{i}}\right)-1,
\end{equation*}
where $x^*=(0,\cdots,0)^T$, $H^*=0$.
\item[(6)] Trigonometric function (n=50, $-50\leq x_i\leq 50$)
\begin{equation*}\label{Tri}
H_6(x)=-\sum_{i=1}^{n}\left[8\sin^2(7(x_i-0.9)^2)+6\sin^2(14(x_i-0.9)^2)+(x_i-0.9)^2\right]-1,
\end{equation*}
where $x^*=(0.9,\cdots,0.9)^T$, $H^*=-1$.
\item[(7)] Rastrigin function (n=20, $-5.12\leq x_i\leq 5.12$)
\begin{equation*}\label{Rastrigin}
H_7(x)=-\sum_{i=1}^{n}\left(x_i^2-10\cos(2\pi x_i)\right)-10n-1,
\end{equation*}
where $x^*=(0,\cdots,0)^T$, $H^*=-1$.
\item[(8)] Pint\'{e}r's function (n=50, $-50\leq x_i\leq 50$)
\begin{eqnarray*}\label{Pinter}
\nonumber
H_8(x)&=&-\left[\sum_{i=1}^{n}ix_i^2+\sum_{i=1}^{n}20i\sin^2(x_{i-1}\sin
x_i-x_i+\sin
x_{i+1})\right.\\
& & \left.
+\sum_{i=1}^{n}i\log_{10}(1+i(x_{i-1}^2-2x_i+3x_{i+1}-\cos
x_i+1)^2)\right]-1,
\end{eqnarray*}
where $x^*=(0,\cdots,0)^T$, $H^*=-1$.
\item[(9)] Levy function (n=50, $-50\leq x_i\leq 50$)
\begin{equation*}\label{Levy}
H_9(x)=-\sin^2(\pi
y_1)-\sum_{i=1}^{n-1}\left[(y_i-1)^2(1+10\sin^2(\pi
y_i+1))\right]-(y_n-1)^2(1+10\sin^2(2\pi y_n))-1,
\end{equation*}
where $y_i=1+ (x_i-1)/4$, $x^*=(1,\cdots,1)^T$, $H^*=-1$.
\item[(10)] Weighted Sphere function (n=50, $-50\leq x_i\leq 50$)
\begin{equation*}\label{Sphere}
H_{10}(x)=-\sum_{i=1}^{n}ix_i^2-1
\end{equation*}
where $x^*=(0,\cdots,0)^T$, $H^*=-1$.
\end{enumerate}
Specifically, Dejong's 5th ($H_1$) and Shekel's ($H_2$) are low-dimensional problems with a small number of local optima that are scattered and far from each other; Powel ($H_3$) and Rosenbrock ($H_4$) are badly-scaled functions; Griewank ($H_5$),
Trigonometric ($H_6$), and Rastrigin ($H_7$) are high-dimensional multimodal problems with a large number of local optima, and the
number of local optima increases exponentially with the problem dimension; Pint\'{e}r ($H_8$) and Levy ($H_9$) are both multimodal
and badly-scaled problems; Weighted Sphere function ($H_{10}$) is a high-dimensional concave function.

We compare the performance of GASS and GASS\_avg with two other algorithms: the modified version of the CE method based on stochastic approximation proposed by \cite{hu:2012a} and the MRAS method proposed by \cite{hu:2007a}. In our comparison, we try to use the same parameter setting in all four methods. The common parameters in all four methods are set as follows: the quantile parameter is set to be $\rho=0.02$ for low-dimensional problems $H_1$ and $H_2$, and $\rho=0.05$ for all the other problems; the parameterized exponential family distribution $f(x;\theta_k)$ is chosen to be independent multivariate normal distribution $\mathcal{N}(\mu_k,\Sigma_k)$; the initial mean $\mu_0$ is chosen randomly according to the uniform distribution on $[-30,30]^n$, and the initial covariance matrix is set to be $\Sigma_0=1000I_{n\times n}$, where $n$ is the dimension of the problem; the sample size at each iteration is set to be $N=1000$. In addition, we observe that the performance of the algorithm is insensitive to the initial candidate solutions if the initial variance is large enough.

In GASS and GASS\_avg, we consider the shape function of the form (\ref{quantile}), i.e.,
\begin{equation*}
S_{\theta_k}(H(x))=(H(x)-H_{lb})\frac{1}{1+e^{-S_0(H(x)-\gamma_{\theta_k})}},
\end{equation*}
In our experiment, $S_0$ is set to be $10^5$, which makes $S_{\theta_k}(H(x))$ a very close approximation to $(H(x)-H_{lb})I\{H(x)\geq \gamma_{\theta_k}\}$; the $(1-\rho)$-quantile
$\gamma_{\theta_k}$ is estimated by the $(1-\rho)$ sample quantile of the function values corresponding to all the candidate solutions generated at the $k^{th}$ iteration. We use the step size: $\alpha_k=\alpha_0/k^{\alpha}$, where $\alpha_0$ reflects the initial step size, and the parameter $\alpha$ should be between $0$ and $1$. We set $\alpha_0=0.3$ for the low-dimensional problems $H_1$ and $H_2$ and the badly-scaled problem $H_4$, and set $\alpha_0=1$ for the rest of the problems; we set $\alpha=0.05$, which is chosen to be relatively small to provide a slowly decaying step size. With the above setting of step size, we can always find a $\beta$ such that the sample size $N_k=1000$ satisfies the Assumption~\ref{assump: gain samplesize}(ii) under a finite number of iterations, e.g. $k<2500$ in our experiment. In GASS\_avg, the feedback weight is $c=0.002$ for problems $H_3$, $H_4$ and $H_8$ and  $c=0.1$ for all other problems.

In the modified CE method, we use the gain sequence $\alpha_k=5/(k+100)^{0.501}$, which is found to work best in the experiments. In the implementation of MRAS method, we use a
smoothing parameter $\nu$ when updating the parameter $\theta_k$ of the parameterized distribution, and set $\nu=0.2$ as suggested by \cite{hu:2007a}. The rest of the parameter setting for MRAS is as follows: $\lambda=0.01$,  $r=10^{-4}$ in the shape function $S(H(x))=\exp\{rH(x)\}$. Other than using an increasing sample size in \cite{hu:2012a} and \cite{hu:2007a}, and updating quantile $\rho_k$
in \cite{hu:2007a}, the constant sample size $N=1000$ and a constant $\rho$ are used in our experiments for a fair comparison of all the methods.

\begin{table*}[ht]
\begin{center}
\begin{minipage}[t]{1\textwidth}
  \centering
\begin{scriptsize}
\begin{tabular}{|c|c|c|c|c|c|c|c|c|c|}
\hline
 & &\multicolumn{2}{c|}{GASS} &\multicolumn{2}{c|}{GASS\_avg}&\multicolumn{2}{c|}{modified CE}&\multicolumn{2}{c|}{MRAS}\\
\cline{2-10}
 &       & \vspace{-0.02in}       &                 & \vspace{-0.02in}       &                 & \vspace{-0.02in} &          & \vspace{-0.02in} &      \\
 & $H^*$ & $\bar{H}^* (std\_err)$ & $M_\varepsilon$ & $\bar{H}^* (std\_err)$ & $M_\varepsilon$ & $\bar{H}^* (std\_err)$ & $M_\varepsilon$& $\bar{H}^* (std\_err)$ & $M_\varepsilon$\\
\hline
Dejong's 5th $H_1$ & -0.998 & -0.998(4.79E-7) & 100 & -0.998(8.97E-7) & 100 & -1.02(0.014) & 95 & -0.9981(6.63E-4) & 98\\
\hline
Shekel $H_2$ & 10.153 & 9.92(0.114) & 96 & 9.91(0.106) & 95 & 10.153(1.09E-7)& 79 & 9.90(0.126) & 96\\
\hline
Powel $H_3$ & -1 & -1(1.48E-6) & 100 & -1(1.89E-6) & 100  & -1(8.87E-9) & 100 & -1.50(0.433) & 95\\
\hline
Rosenbrock $H_4$ & -1 & -1.03(1.40E-4) & 0 & -1.09(0.0301) & 46  & -1.91(0.016) & 0 & -7.10(0.629) & 0\\
\hline
Griewank $H_5$ & 0 & 0(8.45E-15) & 100 & 0(7.30E-15) & 100  & -0(3.02E-16) & 100 & -0.14(0.017) & 57\\
\hline
Trigonometric $H_6$ & -1 & -1(9.72E-13) & 100 & -1(1.08E-12) & 100 & -1(2.23E-18) & 100 & -1(4.69E-7) & 100\\
\hline
Rastrigin $H_7$ & -1 & -1.15(0.0357) & 85 & -1.19(0.044) & 83 & -1.01(0.0099) & 99 & -83.45(0.634) & 0\\
\hline
Pinter $H_8$ & -1 & -1.007(0.0034) & 93 & -1.04(0.0104) & 63  & -6.08(0.0254) & 0 & -530.4(48.64) & 2\\
\hline
Levy $H_9$ & -1 & -1(9.56E-13) & 100 & -1(1.29E-7) & 100 & -1.063(3.87E-18) & 100 & -1(1.42E-10) & 100\\
\hline
Sphere $H_{10}$ & -1 & -1(1.79E-11) & 100 & -1(1.42E-11) & 100 & -1(2.23E-18) & 100 & -1(9.95E-9) & 100\\
\hline
\end{tabular}
\end{scriptsize}
\end{minipage}
\caption{Comparison of GASS, GASS\_avg, modified CE and MRAS}
\label{performance}
\end{center}
\end{table*}

\begin{figure*}[ht]
\begin{center}
\begin{minipage}[t]{.45\textwidth}
  \centering
  {\includegraphics[width=\textwidth]{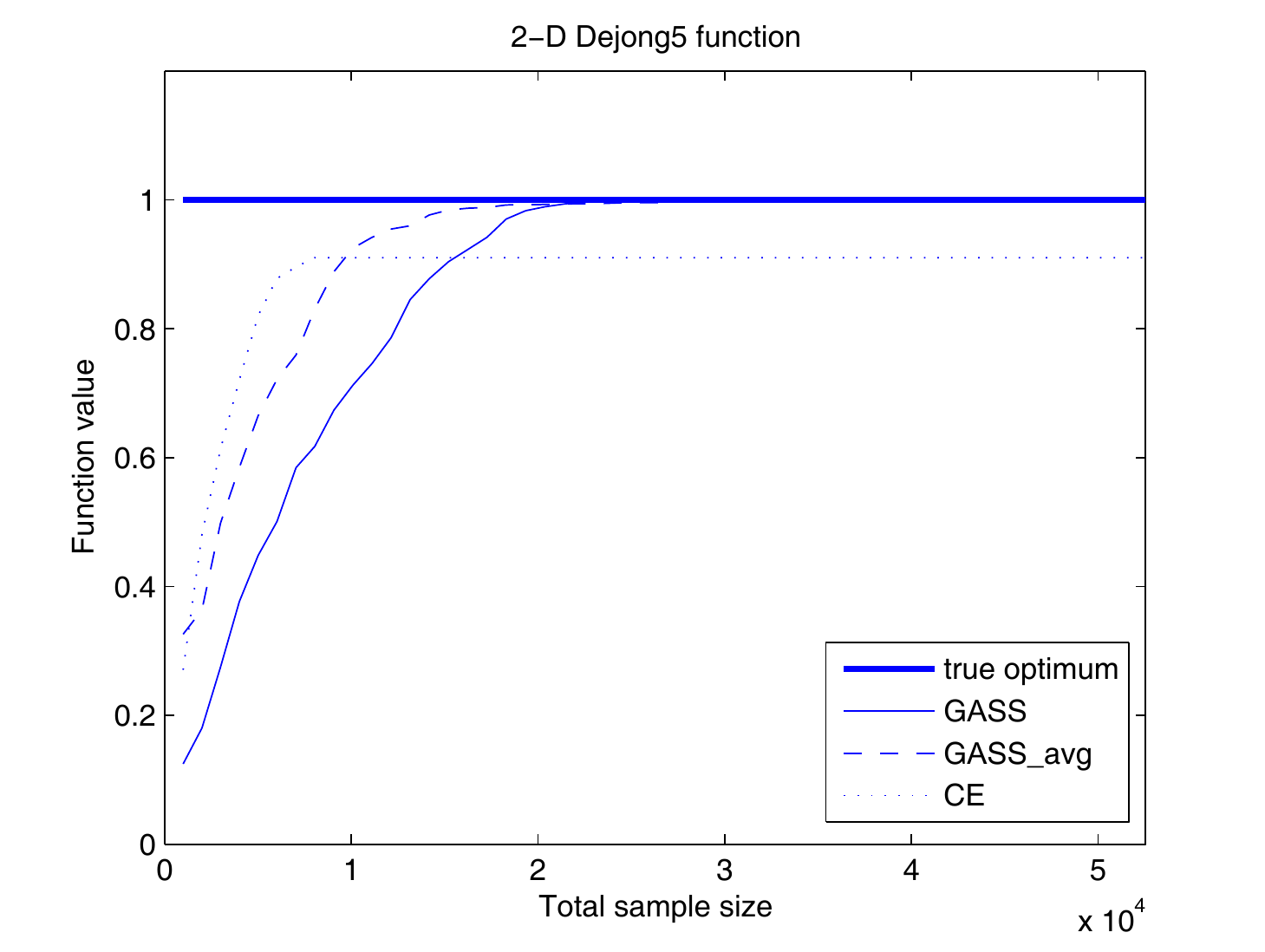}}
\end{minipage}
\begin{minipage}[t]{.45\textwidth}
  \centering
  {\includegraphics[width=\textwidth]{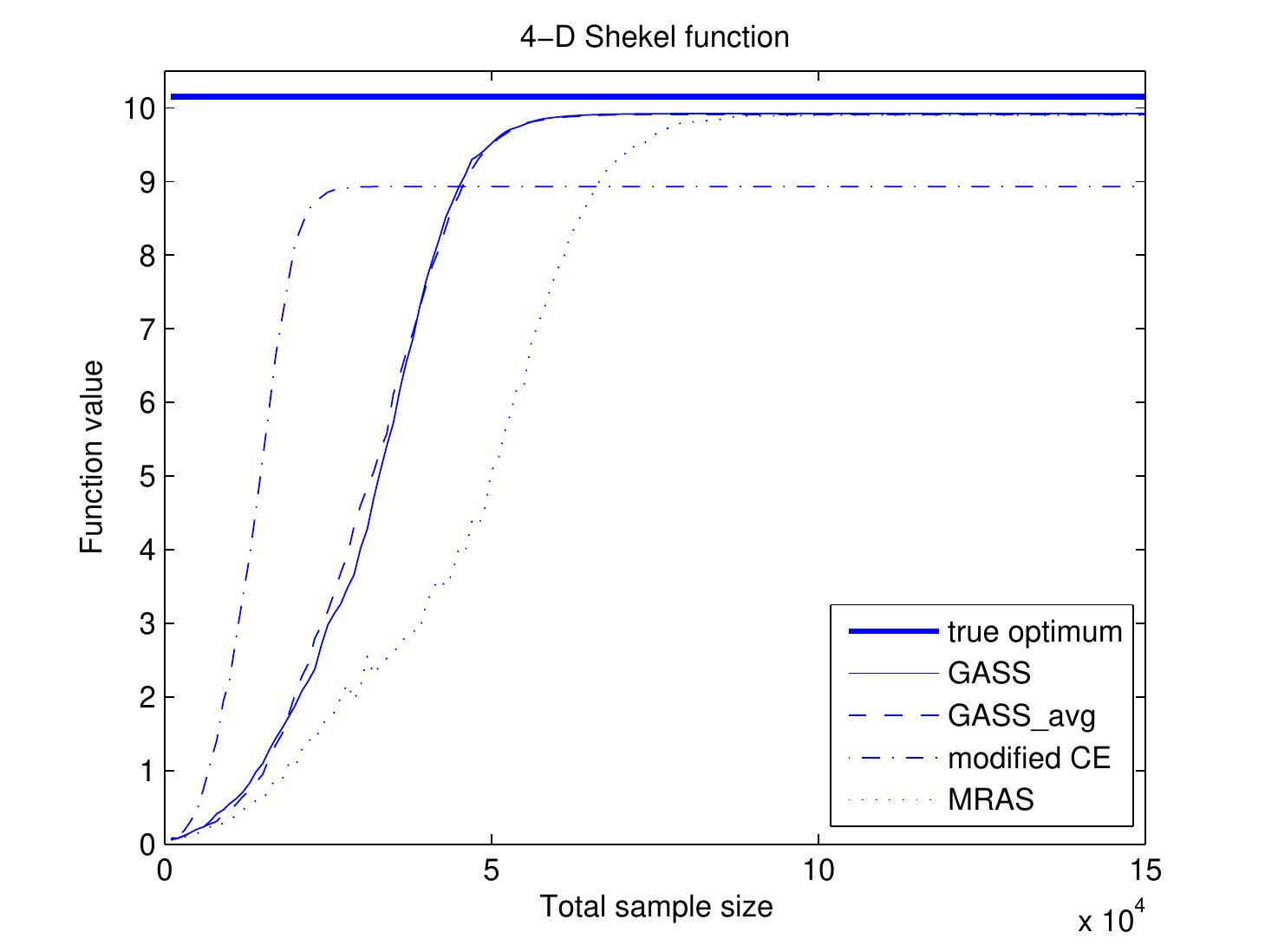}}
\end{minipage}
\begin{minipage}[t]{.45\textwidth}
  \centering
  {\includegraphics[width=\textwidth]{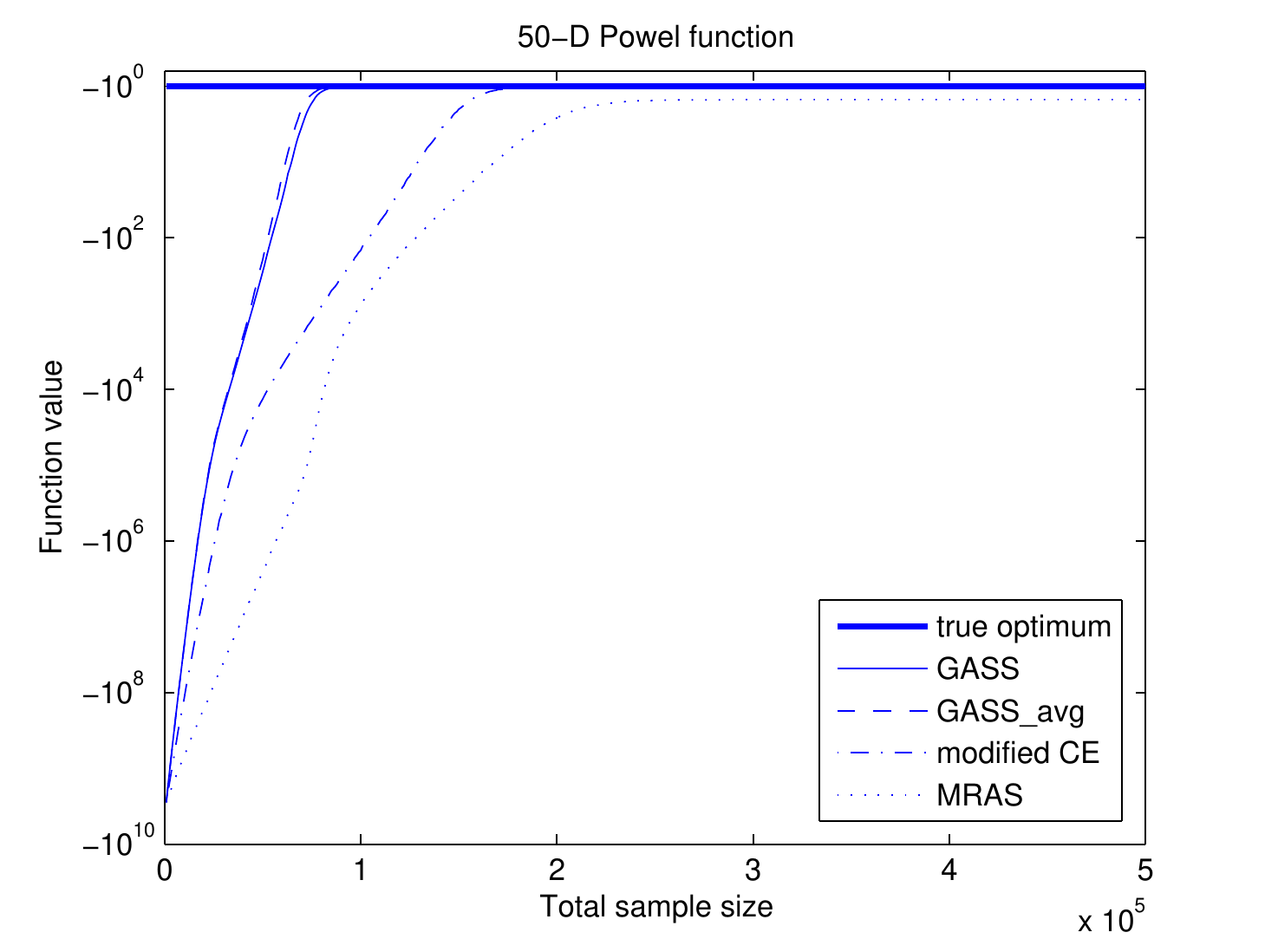}}
\end{minipage}
\begin{minipage}[t]{.45\textwidth}
  \centering
  {\includegraphics[width=\textwidth]{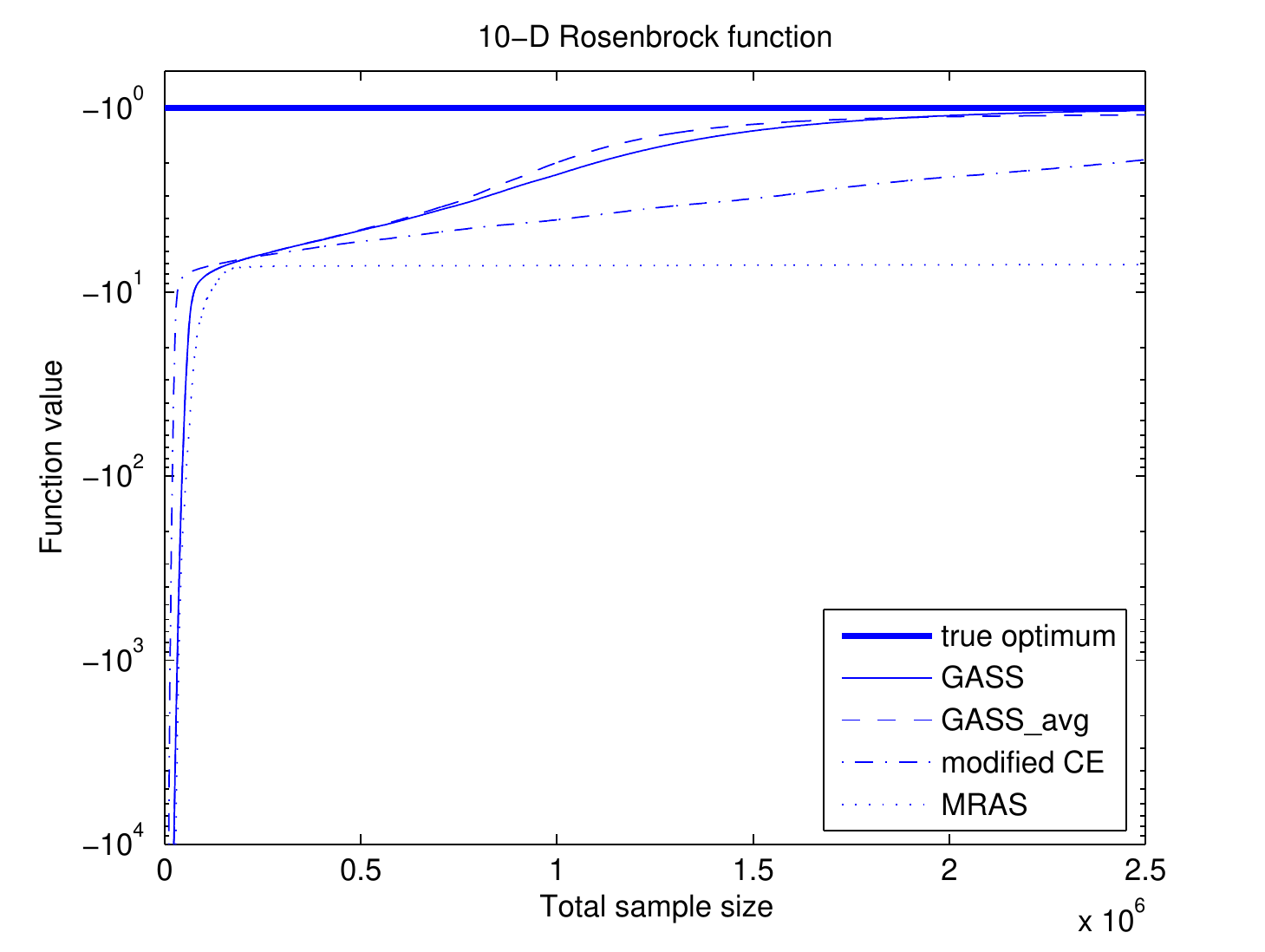}}
\end{minipage}
\begin{minipage}[t]{.45\textwidth}
  \centering
  {\includegraphics[width=\textwidth]{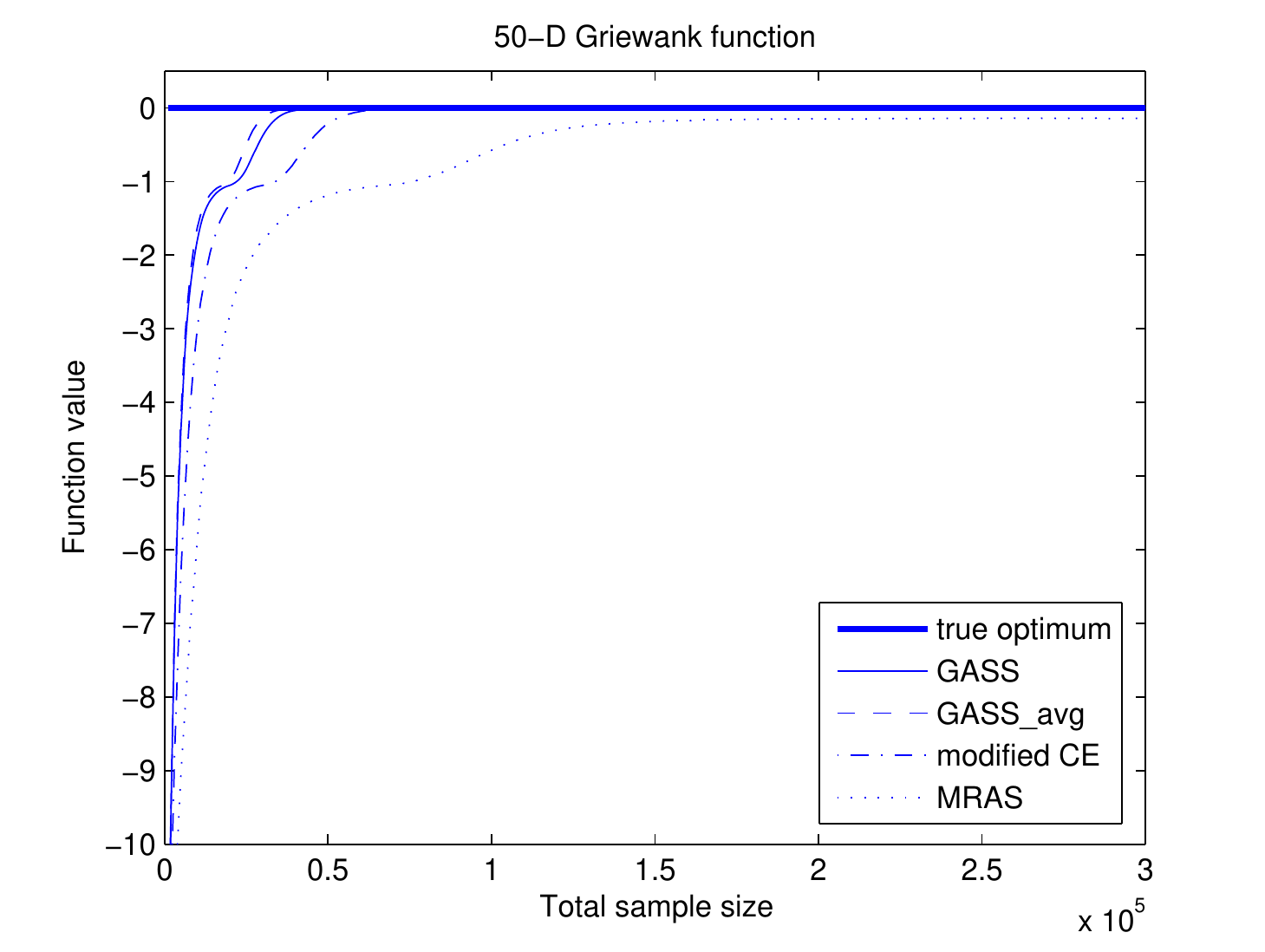}}
\end{minipage}
\begin{minipage}[t]{.45\textwidth}
  \centering
  {\includegraphics[width=\textwidth]{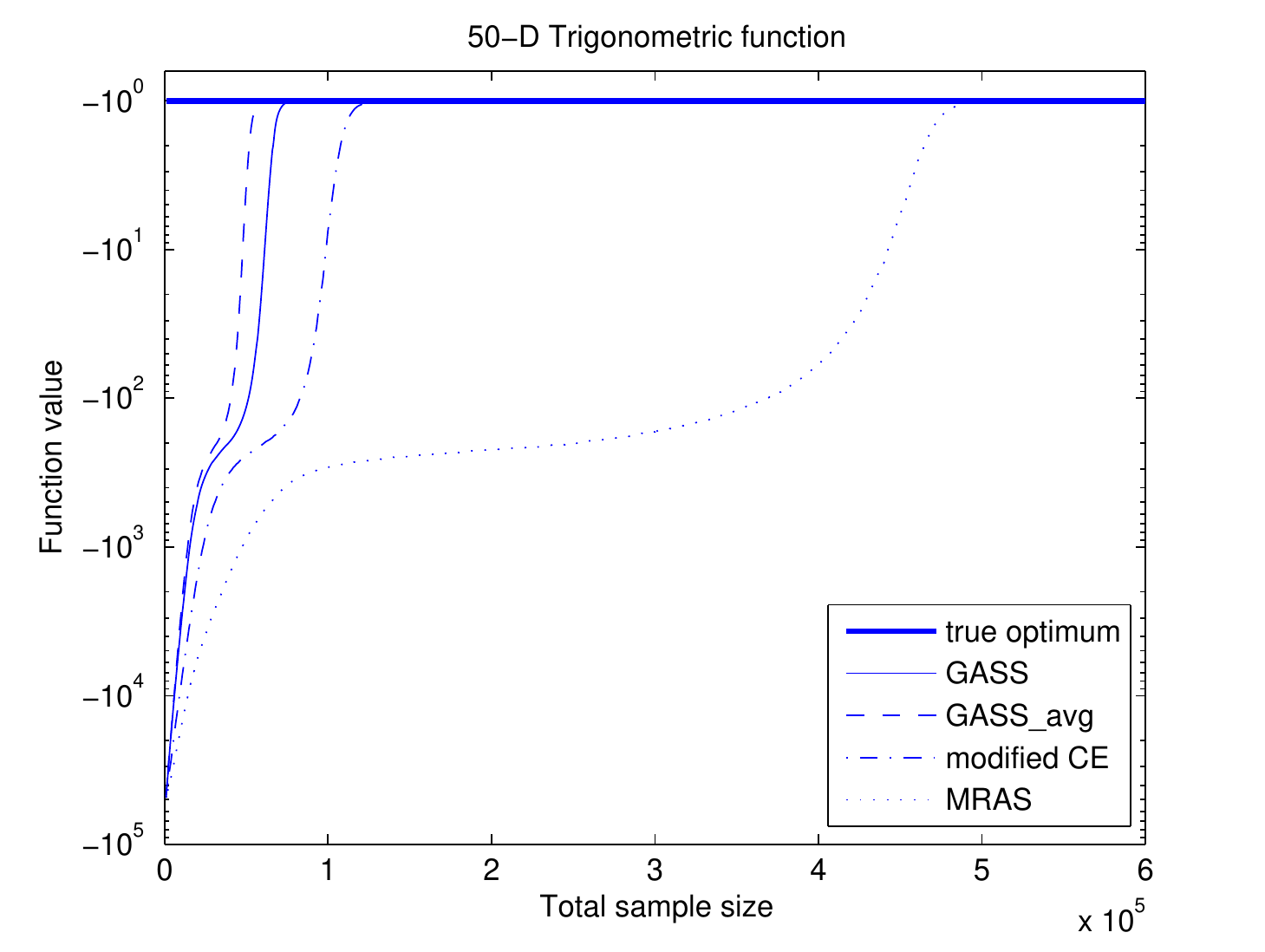}}
\end{minipage}

\caption{Comparison of GASS, GASS\_avg, modified CE and MRAS}
\label{result}
\end{center}
\end{figure*}

\begin{figure*}[ht]
\begin{center}

\begin{minipage}[t]{.45\textwidth}
  \centering
  {\includegraphics[width=\textwidth]{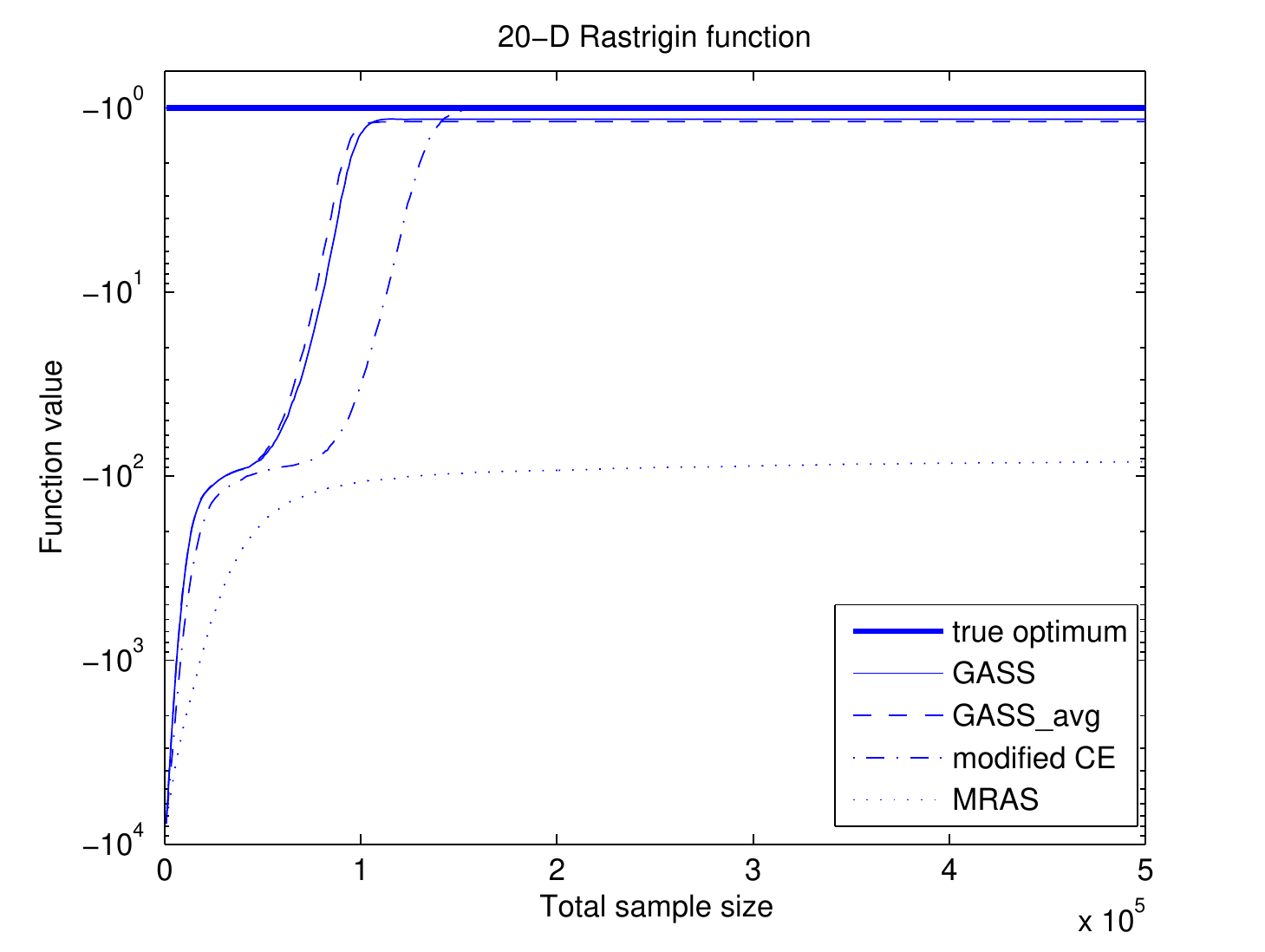}}
\end{minipage}
\begin{minipage}[t]{.45\textwidth}
  \centering
  {\includegraphics[width=\textwidth]{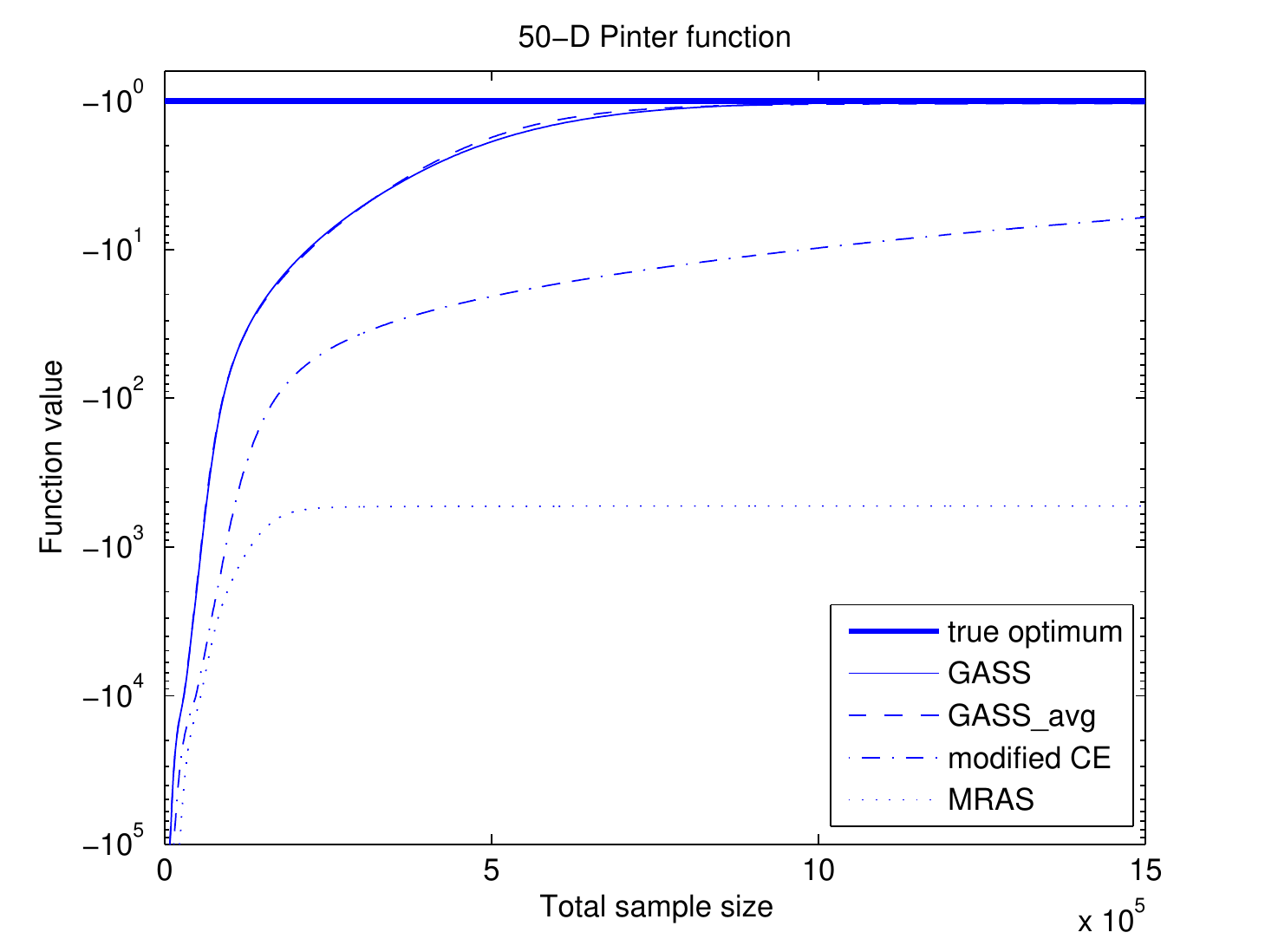}}
\end{minipage}
\begin{minipage}[t]{.45\textwidth}
  \centering
  {\includegraphics[width=\textwidth]{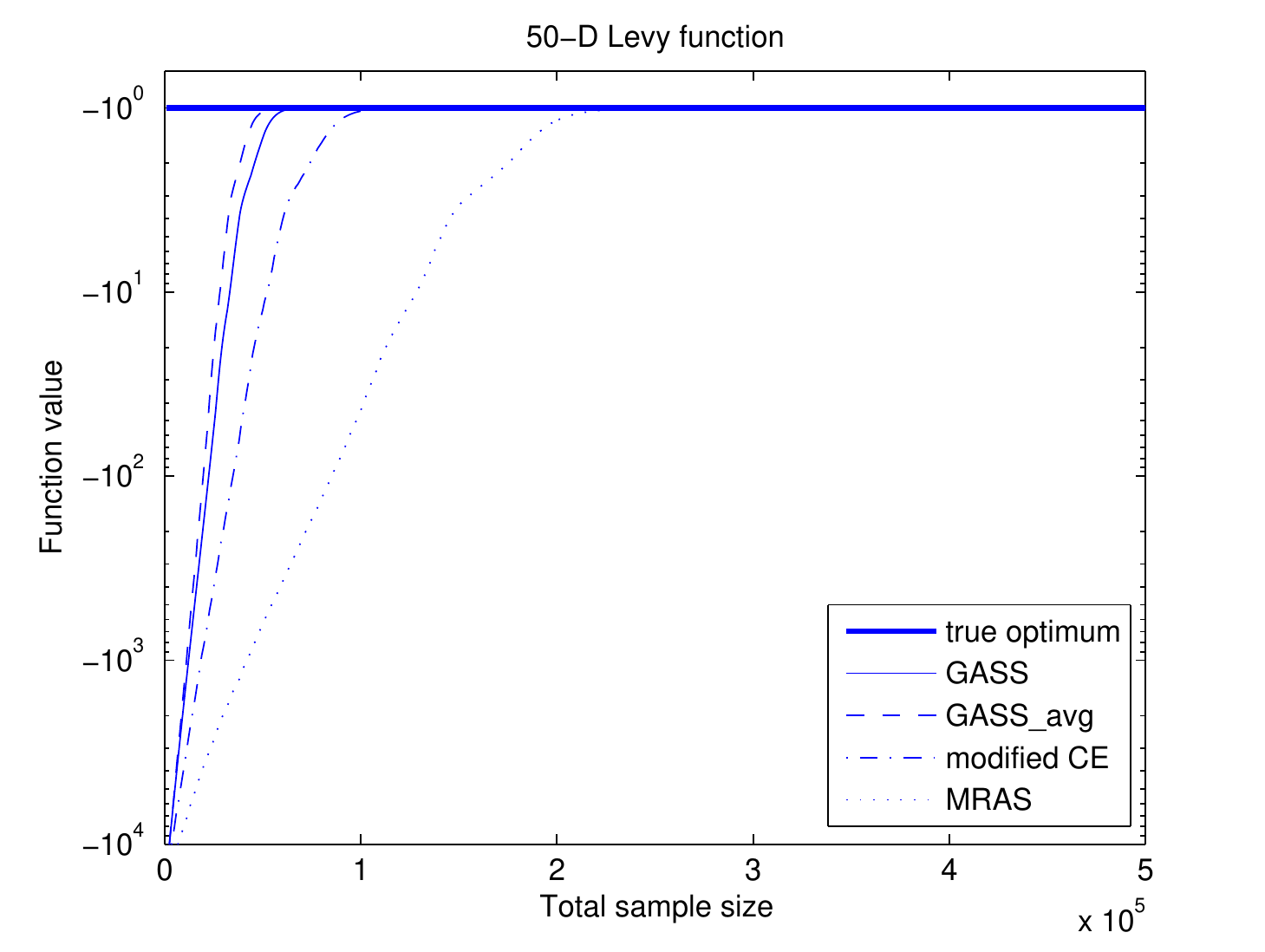}}
\end{minipage}
\begin{minipage}[t]{.45\textwidth}
  \centering
  {\includegraphics[width=\textwidth]{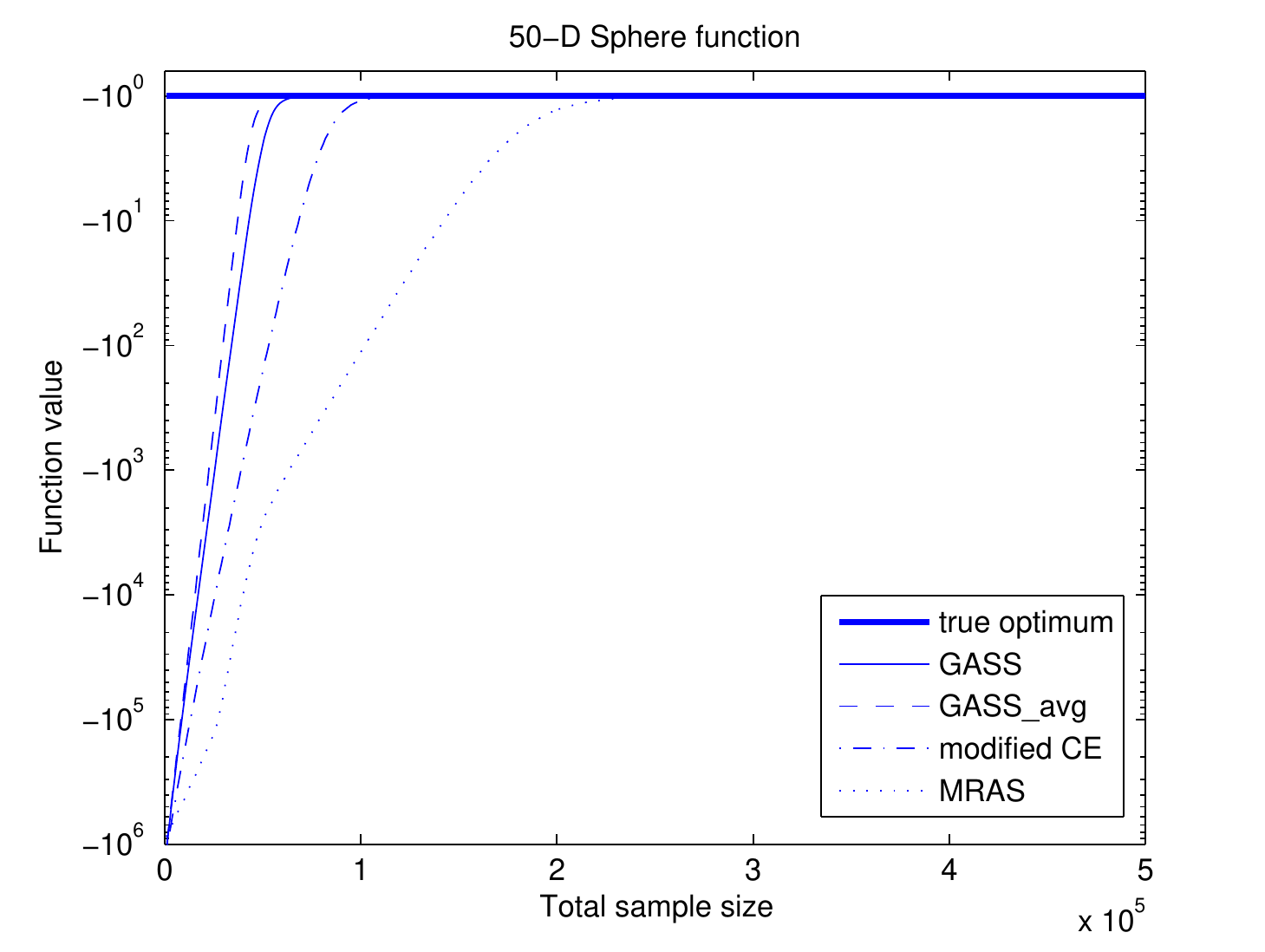}}
\end{minipage}
\caption{Comparison of GASS, GASS\_avg, modified CE and MRAS}
\label{result1}
\end{center}
\end{figure*}

In the experiments, we found the computation time of function
evaluations dominates the time of other steps, so we compare the
performance of the algorithms with respect to the total number of
function evaluations, which is equal to the total number of samples.
The average performance based on 100 independent runs for each
method is shown in Table \ref{performance}, where $H^*$ is the true
optimal value of $H(\cdot)$; $\bar{H}^*$ is the average of the
function values returned by the $100$ runs of an algorithm;
$std\_err$ is the standard error of these $100$  function values;
$M_\varepsilon$ is the number of $\varepsilon$-optimal solutions out
of 100 runs ($\varepsilon$-optimal solution is the solution such
that $H^*-\hat{H}^*\le\varepsilon$, where $\hat{H}^*$ is the optimal
function value returned by an algorithm). We consider
$\varepsilon=10^{-2}$ for problems $H_4$, $H_7$, $H_8$ and
$\varepsilon=10^{-3}$ for all other problems. Fig.~\ref{result} and
Fig.~\ref{result1} show the average (over $100$ runs) of best value of
$H(\cdot)$ at the current iteration versus the total number of samples
generated so far.

From the results, GASS and GASS\_avg find all the $\varepsilon$-optimal solutions in $100$ runs for problems $H_1$,
$H_3$, $H_5$, $H_6$, $H_9$, and $H_{10}$. Modified CE finds all the $\varepsilon$-optimal solutions for problems $H_3$, $H_5$, $H_6$,
$H_9$, and $H_{10}$. MRAS only finds all the $\varepsilon$-optimal solutions for the problems $H_6$ and $H_9$ and the convex problem $H_{10}$.
As for the convergence rate, GASS\_avg always converges faster than GASS, verifying the effectiveness of averaging with online feedback. Both GASS and GASS\_avg converge faster than MRAS on all the
problems, and converge faster than the modified CE method when $\alpha_0$ is set to be large, i.e. on problems $H_3$ and $H_5 - H_{10}$.

\section{Conclusion}\label{sec:6}
In this paper, we have introduced a new model-based stochastic search algorithm for solving general black-box optimization problems. The algorithm generates candidate solutions from a parameterized sampling distribution over the feasible region, and uses a quasi-Newton like iteration on the parameter space of the parameterized distribution to find improved sampling distributions. Thus, the algorithm enjoys the fast convergence speed of classical gradient search methods while retaining the robustness feature of model-based methods. By formulating the algorithm iteration into the form of a generalized stochastic approximation recursion, we have established the convergence and convergence rate results of the algorithm. Our numerical results indicate that the algorithm shows promising performance as compared with some of the existing approaches.

%

\appendix
\section{Appendix}
\noindent
\proof{\textbf{Proposition~\ref{prop: L gradient and hessian}.}}
Consider the gradient of $L(\theta;\theta')$ with respect to $\theta$,
\begin{eqnarray}
\nabla_{\theta}L(\theta; \theta') &=& \int{S_{\theta'}(H(x))\nabla_{\theta}f(x;\theta)dx}  \nonumber \\
&=& \int{S_{\theta'}(H(x))f(x;\theta)\nabla_{\theta}\ln{f(x;\theta)}dx}  \nonumber \\
&=& E_{\theta}[S_{\theta'}(H(X))\nabla_{\theta}\ln{f(X;\theta)}], \label{above11}
\end{eqnarray}
where the interchange of integral and derivative in the first equality follows from the boundedness assumptions on $S_{\theta'}$ and $\nabla_{\theta}f(x;\theta)$ and the dominated convergence theorem.

Consider the Hessian of $L(\theta;\theta')$ with respect to $\theta$,
\begin{eqnarray}
\nabla_{\theta}^2 L(\theta; \theta') &=& \int{S_{\theta'}(H(x))\nabla_{\theta}^2f(x;\theta)dx} \nonumber \\
&=& \int{S_{\theta'}(H(x))f(x;\theta)\nabla_{\theta}^2\ln{f(x;\theta)}dx} + \int{S_{\theta'}(H(x))\nabla_{\theta}\ln{f(x;\theta)}\nabla_{\theta}{f(x;\theta)^T}dx} \nonumber\\
&=& E_{\theta}[S_{\theta'}(H(X))\nabla_{\theta}^2\ln{f(X;\theta)}] + E_{\theta}[S_{\theta'}(H(X))\nabla_{\theta}\ln{f(x;\theta)}\nabla_{\theta}\ln{f(x;\theta)}^T], \label{above13}
\end{eqnarray}
where the last equality follows from the fact that $\nabla_{\theta}f(x;\theta) = f(x;\theta)\nabla_{\theta}\ln f(x;\theta)$.

Furthermore, if $f(x;\theta) = \exp\{\theta^T T(x) - \phi(\theta) \}$, we have
\begin{eqnarray}
\nabla_{\theta} \ln{f(x;\theta)} &=& \nabla_{\theta}\left(\theta^T T(x) -  \ln\int{\exp(\theta^TT(x))dx} \right) \nonumber \\
&=& T(x) - \frac{\int{\exp(\theta^T T(x)) T(x)dx}} {\int{\exp(\theta^T T(x))dx}} \nonumber \\
&=& T(x) - E_{\theta}[T(X)]. \label{above2}
\end{eqnarray}
Plugging (\ref{above2}) into (\ref{above11}) yields
$$
\nabla_{\theta} L(\theta; \theta')  = E_{\theta}[S_{\theta'}(H(X))T(X)] - E_{\theta}[S_{\theta'}(H(X))] E_{\theta}[T(X)].
$$

Differentiating (\ref{above2})  with respect to $\theta$, we obtain
\begin{eqnarray}
\nabla_{\theta}^2 \ln{f(x;\theta)} &=& - \frac{\int \exp(\theta^T T(x))T(x)T(x)^Tdx}{\int{\exp(\theta^T T(x))dx}} \nonumber \\
&& +~  \frac{\int \exp(\theta^T T(x))T(x)dx \left( \int \exp(\theta^T T(x))T(x)dx \right)^T }{ \left( \int{\exp(\theta^T T(x))dx} \right)^2}  \nonumber \\
&=& -E_{\theta}[T(X)T(X)^T] + E_{\theta}[T(X)] E_{\theta}[T(X)]^T \nonumber \\
&=& - \mathrm{Var}_{\theta}[T(X)]. \label{above4}
\end{eqnarray}
Plugging (\ref{above2}) and (\ref{above4}) into (\ref{above13}) yields
\begin{eqnarray*}
\nabla_{\theta}^2{L(\theta; \theta')} &=& E_{\theta}[S_{\theta'}(H(X))(T(X)-E_{\theta}[T(X)])(T(X)-E_{\theta}[T(X)])^T]  \nonumber \\
  && -~ \mathrm{Var}_{ \theta}[T(X)]E_{\theta}[S_{\theta'}(H(X))].
\end{eqnarray*}
\endproof

\proof{\textbf{Proposition~\ref{prop: gradient and hessian}.}}
Consider the gradient of $l(\theta; \theta')$ with respect to $\theta$,
\begin{eqnarray}
\nabla_{\theta}l(\theta; \theta')|_{\theta = \theta'} &=&  \frac{\nabla_{\theta}L(\theta; \theta')}{L(\theta;\theta')} \bigg|_{\theta = \theta'} \nonumber \\
&=&  \frac{\int{ S_{\theta'}(H(x)) f(x;\theta) \nabla_{\theta}\ln{f(x; \theta)}dx}} {L(\theta;\theta')}\bigg|_{\theta = \theta'} \label{above1}  \\
&=& E_{p(\cdot; \theta')} [\nabla_{\theta}\ln{f(X;\theta')}]. \nonumber
\end{eqnarray}

Differentiating (\ref{above1}) with respect to $\theta$, we obtain the Hessian
\begin{eqnarray*}
\nabla_{\theta}^2 l(\theta;\theta')|_{\theta = \theta'} &=& \frac{\int{S_{\theta'}(H(x))f(x;\theta)\nabla_{\theta}^2 \ln{f(x;\theta)}dx}}{L(\theta;\theta')}  + \frac{\int{S_{\theta'}(H(x))\nabla_{\theta}\ln f(x;\theta)(\nabla_{\theta}f(x;\theta))^Tdx}}{L(\theta; \theta')}...\\
&& -~ \frac{(\int{S_{\theta'}(H(x))f(x;\theta)\nabla_{\theta}\ln f(x;\theta)dx}) (\nabla_\theta L(\theta;\theta'))^T}{L(\theta;\theta')^2}\bigg|_{\theta = \theta'}
\end{eqnarray*}
Using $\nabla_{\theta}f(x;\theta) = f(x;\theta)\nabla_{\theta}\ln f(x;\theta)$ in the second term on the righthand side, the above expression can be written as
\begin{eqnarray}
\nabla_{\theta}^2 l(\theta; \theta')|_{\theta = \theta'} &=& E_{p(\cdot; \theta')}[\nabla_{\theta}^2 \ln{f(X;\theta')}]  + E_{p(\cdot ; \theta')}\left[ \nabla_{\theta'}\ln{f(X;\theta')} (\nabla_{\theta'}\ln{f(X;\theta')})^T \right]  \nonumber \\
&& - ~E_{p(\cdot ; \theta')}\left[ \nabla_{\theta}\ln{f(X;\theta')}\right] E_{p(\cdot ; \theta')}\left[ \nabla_{\theta}\ln{f(X;\theta')}\right] ^T  \nonumber \\
&=&  E_{p(\cdot; \theta')}[\nabla_{\theta}^2 \ln{f(X;\theta')}]  + \mathrm{Var}_{p(\cdot ; \theta')}\left[ \nabla_{\theta}\ln{f(X;\theta')} \right]. \label{above3}
\end{eqnarray}

Furthermore, if $f(x;\theta) = \exp\{\theta^T T(x) - \phi(\theta) \}$,  plugging (\ref{above2}) into (\ref{above1}) yields
$$
\nabla_{\theta} l(\theta; \theta')|_{\theta = \theta'}  = E_{p(\cdot; \theta')}[T(X)] - E_{\theta'}[T(X)],
$$
and plugging (\ref{above2}) and (\ref{above4}) into (\ref{above3}) yields
$$
\nabla_{\theta}^2{l(\theta;\theta')}|_{\theta = \theta'} = \mathrm{Var}_{p(\cdot; \theta')}[T(X)] - \mathrm{Var}_{\theta'}[T(X)].
$$
\endproof

\noindent \proof{\textbf{Lemma~\ref{lem:quant}.}} Because $S_{\theta}$ is continuous in $\gamma_{\theta}$, it is sufficient to show that $\widehat{\gamma}_{\theta_k} \rightarrow \gamma_{\theta_k}$ w.p.1 as $k \rightarrow \infty$, which can be shown in the same way as Lemma~7 in \cite{hu:2007a}, except that we need to verify the following condition in their proof:
$$
\sum_{k=1}^{\infty}\exp\left(-\tilde{M}N_k \right) < \infty,
$$
where $\tilde{M}$ is positive constant. It is easy to see that this condition is
trivially satisfied in our setting by taking \ez{$N_k=N_0k^{\zeta}$} with $\zeta>0$.
\endproof

\noindent \proof{\textbf{Lemma~\ref{lem:1}.}}
Before the formal proof of Lemma~\ref{lem:1}, we first introduce a key inequality to our proof - the matrix bounded differences inequality (\cite{tropp:2011}), which is a matrix version of the generalized Hoeffding inequality (i.e., McDiarmid's inequality (\cite{mcdiarmid:1989})). Let $\lambda_{\max}(\cdot)$ and $\lambda_{\min}(\cdot)$ return the largest and smallest eigenvalue of a matrix, respectively.
\begin{Theorem}(Matrix bounded differences, Corollary~7.5, \cite{tropp:2011}) \label{thm:matrix}
Let $\{X^i: i = 1,2,\ldots,N\}$ be an independent family of random variables, and let $V$ be a function that maps $N$ variables to a self-adjoint matrix of dimension $d$. Consider a sequence of $\{C_k\}$ of fixed self-adjoint matrices that satisfy
$$
\left( V(x^1,\ldots, x^i,\ldots,x^N) - V(x^1,\ldots, \tilde{x}^i,\ldots,x^N) \right)^2 \leq C_i^2,
$$
where $x^i$ and $\tilde{x}^i$ range over all possible values of $X^i$ for each index $i$. Compute the variance parameter
$$
\sigma^2 := \left\|\sum_i C_i^2\right\|_2.
$$
Then, for all $\delta > 0$,
$$
P\left\{\lambda_{\max}(V(\mathbf{x}) - E[V(\mathbf{x})])\geq \delta \right\} \leq d\exp\left\{\frac{-\delta^2}{8\sigma^2}\right\},
$$
where $\mathbf{x} = (X^1,\ldots,X^N).$
\end{Theorem}

Now we proceed to the formal proof of Lemma~\ref{lem:1}. Recall that $\xi_k$ can be written as
\begin{equation} \label{xi}
\xi_k = (\widehat{V}_k^{-1} - V_k^{-1})(\widetilde{E}_{p_k}[T(X)]-E_{\theta_k}[T(X)]) + V_k^{-1}(\widetilde{E}_{p_k}[T(X)] - E_{p_k}[T(X)]).
\end{equation}
To bound the first term on the right-hand-side in (\ref{xi}), we notice that since $V_k^{-1}$ and $\widehat{V}_k^{-1}$ are both positive definite and $(\epsilon^{-1}I - V_k^{-1})$ and $(\epsilon^{-1}I - \widehat{V}_k^{-1})$ are both positive semi-definite, we have
\begin{eqnarray}
\|V_k^{-1}-\widehat{V_k}^{-1}\| &=& \|V_k^{-1}(\widehat{V}_k - V_k)\widehat{V}_k^{-1}\| \nonumber \\
&\leq&  \|V_k^{-1}\| \|\widehat{V}_k - V_k\| \|\widehat{V}_k^{-1}\| \nonumber \\
&\leq& \epsilon^{-2} \| \widehat{V}_k - V_k\|.  \label{above7}
\end{eqnarray}
To establish a bound on $\|\widehat{V}_k- V_k\|$, we use the matrix bounded differences inequality that is introduced above. For simplicity of exposition, we drop the subscript $k$ in the expression below.
\begin{align*}
\begin{split}
& \sup_{x^i, \tilde{x}^i \in \mathcal{X}}\left\{\widehat{V}(x^1,\ldots, x^i, \ldots, x^N) - \widehat{V}(x^1,\ldots, \tilde{x}^i, \ldots, x^N)\right\}^2 \\
&= \frac{1}{N^2}\sup_{x^i, \tilde{x}^i \in \mathcal{X}}\left\{\left[T(x^i)T(x^i)^T - T(\tilde{x}^i)T(\tilde{x}^i)^T\right] - \frac{1}{N-1}\sum_{j\neq i}{\left(T(x^i)-T(\tilde{x}^i)\right)T(x^j)^T} ...
 \parenthnewln - \frac{1}{N-1}\sum_{j\neq i}{T(x^j)\left(T(x^i)-T(\tilde{x}^i)\right)^T} \right\}^2 \\
&\leq \frac{1}{N^2}C,
\end{split}
\end{align*}
where $C$ is a fixed positive semidefinite matrix. This last inequality is due to Assumption~\ref{assump: gain samplesize}(iv) that $T(x)$ is bounded on $\mathcal{X}$.  Note that conditioning on $\mathcal{F}_{k-1}$, $\{x_k^i, i=1,\ldots,N_k\}$ are i.i.d., and $E_{\theta_k}[\widehat{V}_k|\mathcal{F}_{k-1}] = V_k$. Then according to the matrix bounded differences inequality, for all $\delta > 0$,
$$
P\left\{\lambda_{\max}(\widehat{V}_k - V_k) \geq \delta\ | \mathcal{F}_{k-1}\right\} \leq d \exp{\left(\frac{-N_k \delta^2}{8\|C\|_2} \right)},
$$
which also implies
$$
P\left\{-\lambda_{\min}(\widehat{V}_k - V_k) \geq \delta\ | \mathcal{F}_{k-1}\right\} = P\left\{\lambda_{\max}(V_k - \widehat{V}_k) \geq \delta\ | \mathcal{F}_{k-1}\right\} \leq d \exp{\left(\frac{-N_k \delta^2}{8\|C\|_2} \right)}.
$$
Recall that for a symmetric matrix $A$, $\|A\|_2  = \max(\lambda_{\max}(A), -\lambda_{\min}(A))$ and $\|A\| \leq \|A\|_2$. Hence,
$$
P\left\{\|\widehat{V}_k - V_k\| \geq \delta\ | \mathcal{F}_{k-1} \right\} \leq P\left\{\|\widehat{V}_k - V_k\|_2 \geq \delta \ | \mathcal{F}_{k-1} \right\} \leq 2d \exp{\left(\frac{-N_k \delta^2}{8\|C\|_2} \right)}.
$$
Recall that for any nonnegative random variable $X$,
\begin{eqnarray*}
E[X] &=& \int_{0}^{\infty}{P(X \geq x)dx} \\
&\leq& a + \int_{a}^{\infty}{P\left(X \geq x \right)dx}.
\end{eqnarray*}
So we have
\begin{eqnarray*}
E\left[\|\widehat{V}_k - V_k\|^2\ |\mathcal{F}_{k-1}\right] &\leq& a + \int_{a}^{\infty}{P\left\{\|\widehat{V}_k - V_k\|\geq \sqrt{x} \ | \mathcal{F}_k \right\}dx} \\
&\leq& a + \int_{a}^{\infty}{2d \exp{\left(\frac{-N_k x}{8\|C\|_2} \right)}dx}.
\end{eqnarray*}
Set $a = 8\|C\|_2\log{(2d)}/N_k$, and we obtain
\begin{equation} \label{above5}
E\left[\|\widehat{V}_k - V_k\|\ |\mathcal{F}_{k-1}\right]^2 \leq E\left[\|\widehat{V}_k - V_k\|^2\ |\mathcal{F}_{k-1}\right] \leq \frac{8\|C\|_2(1 + \log{(2d)})}{N_k}.
\end{equation}

To bound the second term in the right-hand-side of (\ref{xi}), notice that $\widetilde{E}_{p_k}[T_j(X)]$ is a self-normalized importance sampling estimator of $E_{p_k}[T_j(X)]$, where $T_j(X)$ is the $j^{th}$ element in the vector $T(X)$. Applying Theorem 9.1.10 (pp. 294, \cite{cappe:2005}), we have
$$
E\left[|\widetilde{E}_{p_k}[T_j(X)] - E_{p_k}[T_j(X)]|^2|\mathcal{F}_{k-1}\right] \leq \frac{c_j}{N_k}, ~~j=1,\ldots,d,
$$
where $c_j$'s are positive constants due to the boundedness of $T_j(x)$ on $\mathcal{X}$. Hence,
\begin{eqnarray}
&& E\left[\|\widetilde{E}_{p_k}[T(X)] - E_{p_k}[T(X)]\| |\mathcal{F}_{k-1}\right]^2 \nonumber \\
&\leq& E\left[\|\widetilde{E}_{p_k}[T(X)] - E_{p_k}[T(X)]\|^2 |\mathcal{F}_{k-1}\right] \nonumber \\
&\leq& \sum_{j=1}^{d} E\left[|\widetilde{E}_{p_k}[T_j(X)] - E_{p_k}[T_j(X)]|^2|\mathcal{F}_{k-1}\right] \leq \frac{d\max_jc_j}{N_k}. \label{above6}
\end{eqnarray}

Putting (\ref{above5}) and (\ref{above6}) together, we obtain
\begin{eqnarray*}
E[\|\xi_k\|] &\leq& E\left[\epsilon^{-2}\|\widehat{V}_k-V_k\|\|\widetilde{E}_{p_k}[T(X)]-E_{\theta_k}[T(X)]\| + \|V_k^{-1}\|\|\widetilde{E}_{p_k}[T(X)]-E_{p_k}[T(X)]\| \right] \\
&\leq& M\epsilon^{-2}E\left[E\left[\|\widehat{V}_k-V_k\||\mathcal{F}_{k-1}\right]\right] + \epsilon^{-1}E\left[E\left[\|\widetilde{E}_{p_k}[T(X)]-E_{p_k}[T(X)]\||\mathcal{F}_{k-1}\right]\right]\\
&\leq& \frac{M\epsilon^{-2}\sqrt{8\|C\|_2(1 + \log{(2d)})}+ \epsilon^{-1}\sqrt{d\max_jc_j}}{\sqrt{N_k}} \\
&\triangleq& \frac{c}{\sqrt{N_k}},
\end{eqnarray*}
where the positive constant $M$ is due to the boundedness of $T(x)$ on $\mathcal{X}$.

Therefore, for any $T>0$
\begin{eqnarray*}
E\left[\sum_{i=k}^{\infty}{\alpha_i \|\xi_i\|}\right] &=& \sum_{i=k}^{\infty}{\alpha_i E[\|\xi_i\|]} \\
&\leq& c \sum_{i=k}^{\infty}{\frac{\alpha_i}{\sqrt{N_i}}} \\
&=& c \sum_{i=k}^{\infty}{\frac{1}{i^{\beta}}}  \\
&\leq& c \left( \frac{1}{k^{\beta}} + \int_{k}^{\infty}{\frac{1}{x^{\beta}}dx} \right) \\
&=& c \left(\frac{1}{k^{\beta}} + \frac{1}{\beta -1} \frac{1}{k^{\beta-1}}\right),
\end{eqnarray*}
where the first line follows from the monotone convergence theorem, and the third line follows from Assumption~\ref{assump: gain samplesize}(ii). For any $\tau > 0$, we have from Markov's inequality
\begin{eqnarray*}
P\left(\sum_{i=k}^{\infty}{\alpha_i\|\xi_i\|} \geq \tau\right) &\leq& \frac{E\left[\sum_{i=k}^{\infty}{\alpha_i\|\xi_i\|}\right]}{\tau} \\
&\leq& \frac{c}{\tau }\left(\frac{1 }{k^{\beta}} + \frac{1}{\beta-1}\frac{1}{k^{\beta-1}} \right) ~\rightarrow~ 0 ~~\text{as}~k\rightarrow \infty,
\end{eqnarray*}
where the last statement is due to $\beta>1$. This result of convergence in probability together with the fact that the sequence $\{\sum_{i=k}^{\infty}{\alpha_i\|\xi_i\|}\}$ is monotone implies that the sequence $\{\sum_{i=k}^{\infty}{\alpha_i\|\xi_i\|}\}$ converges w.p.1 as $k \rightarrow \infty$. Furthermore, since $\sup_{\{n:0\leq \sum_{i=k}^{n-1}\alpha_i\leq T\}}\|\sum_{i=k}^{n}{\alpha_i\xi_i}\| \leq \sup_{\{n:0\leq \sum_{i=k}^{n-1}\alpha_i\leq T\}}\sum_{i=k}^{n}{\alpha_i\|\xi_i\|} \leq \sum_{i=k}^{\infty}{\alpha_i\|\xi_i\|}$, we conclude that $\{\sup_{\{n:0\leq \sum_{i=k}^{n-1}\alpha_i\leq T\}}\|\sum_{i=k}^{n}{\alpha_i\xi_i}\|\}$ converges to $0$ w.p.1 as $k \rightarrow \infty$.
\endproof

\proof{\textbf{Theorem~\ref{thm:SA}.}}
To show our theorem, we apply \ez{Theorem~2.1 in \cite{kushner:2010}. The condition on the step size sequence in their theorem is satisfied by our Assumption~\ref{assump: gain samplesize}(i), and condition (2.2) there is a result of Lemma~\ref{lem:1}. Thus, to establish convergence,} it is sufficient to show $b_k \rightarrow 0$ w.p.1 as $k\rightarrow \infty$. Note that
\begin{eqnarray*}
b_k &=& \widehat{V}_k^{-1}\left(\widehat{E}_{p_k}[T(X)] - \widetilde{E}_{p_k}[T(X)]\right) \\
&=& \widehat{V}_k^{-1} \left(\frac{\bar {\mathbb{U}}_k}{\bar {\mathbb{V}}_k}-\frac{\bar {\mathbb{U}}_k}{\tilde {\mathbb{V}}_k}+\frac{\bar {\mathbb{U}}_k}{\tilde {\mathbb{V}}_k}-\frac{\tilde {\mathbb{U}}_k}{\tilde {\mathbb{V}}_k} \right) \\
&=& \widehat{V}_k^{-1} \bar {\mathbb{U}}_k \Big(\frac{\tilde {\mathbb{V}}_k-\bar {\mathbb{V}}_k}{\bar {\mathbb{V}}_k \tilde {\mathbb{V}}_k} \Big)+ \widehat{V}_k^{-1}\frac{\bar {\mathbb{U}}_k-\tilde {\mathbb{U}}_k}{\tilde {\mathbb{V}}_k}.
\end{eqnarray*}
Hence,
\begin{eqnarray*}
\|b_k\| &\leq& \frac{\|\widehat{V}_k^{-1}\|\|\bar {\mathbb{U}}_k \|}{|\bar {\mathbb{V}}_k \tilde {\mathbb{V}}_k|}|\tilde {\mathbb{V}}_k-\bar {\mathbb{V}}_k |+ \frac{\|\widehat{V}_k^{-1}\|}{|\tilde {\mathbb{V}}_k|}\|\bar {\mathbb{U}}_k-\tilde {\mathbb{U}}_k\|\\ \nonumber
&\leq& \frac{\|\widehat{V}_k^{-1}\|\|\bar {\mathbb{U}}_k \|}{|\bar {\mathbb{V}}_k \tilde {\mathbb{V}}_k|}\frac{1}{N_k}\sum_{i=1}^{N_k}|\widehat S_{\theta_k}(H(x^i_k))-S_{\theta_k}(H(x^i_k)) |\\
&& +\frac{\|\widehat{V}_k^{-1}\|}{|\tilde {\mathbb{V}}_k|}\frac{1}{N_k}\sum_{i=1}^{N_k}|\widehat S_{\theta_k}(H(x^i_k))-S_{\theta_k}(H(x^i_k)) |\|T(x^i_k)\|.
\end{eqnarray*}
Since $T(x)$ is bounded, it is easy to see that $\frac{\|\bar {\mathbb{U}}_k \|}{|\bar {\mathbb{V}}_k |}$ is also bounded. Furthermore, note that $\|\widehat{V}_k^{-1}\|$ is bounded and  $|\tilde {\mathbb{V}}_k|$ is bounded away from zero. This together with Assumption~\ref{assump: gain samplesize}(iv) imply that the sequence $\{b_k\}$ converges to zero w.p.1.
\endproof

\noindent \proof{\textbf{Lemma~\ref{lem1}.}}
Under Assumption 1, we know that the sequence $\{\theta_k\}$ converges w.p.1. to a limiting point $\theta^*$. This, together with Assumption~\ref{assump: gain samplesize}(iii), implies that the sequence of sampling distributions $\{f(x;\theta_k)\}$ will converge point-wise in $x$ to a limiting distribution $f(x;\theta^*)$ w.p.1. Note that
$\|\widehat {\mathrm{Var}}_{\theta_k}(T(X))-\mathrm{Var}_{\theta^*}(T(X))\|\leq \|\widehat {\mathrm{Var}}_{\theta_k}(T(X))-\mathrm{Var}_{\theta_k}(T(X)) \| + \|\mathrm{Var}_{\theta_k}(T(X))-\mathrm{Var}_{\theta^*}(T(X)) \|$. Clearly, the first term converges to zero by the strong consistency of the variance estimator. On the other hand, using the point-wise convergence of $\{f(\cdot;\theta_k)\}$ and the dominated convergence theorem, it is easy to see that the second term also vanishes to zero. This shows $\Phi_k\rightarrow \Phi$ w.p.1. Thus, the convergence of $\Gamma_k$ to $\Gamma$ is a direct consequence of the continuity assumption of $J_\mathcal{L}$ in the neighborhood of $\theta^*$. Regarding $T_k$, we have
\begin{align*}
T_k & =  k^{\tau/2} \Phi_k \Big(\frac{\bar {\mathbb{U}}_k}{\bar {\mathbb{V}}_k}-\frac{\bar {\mathbb{U}}_k}{\tilde {\mathbb{V}}_k}+\frac{\bar {\mathbb{U}}_k}{\tilde {\mathbb{V}}_k}-\frac{\tilde {\mathbb{U}}_k}{\tilde {\mathbb{V}}_k} \Big)+k^{\tau/2}\Phi_k\Big(E_{\theta_k}\Big[\frac{\tilde {\mathbb{U}}_k}{\tilde {\mathbb{V}}_k} \Big|\mathcal{F}_{k-1}\Big]-\frac{{\mathbb{U}}_k}{{\mathbb{V}}_k} \Big)\\
&=T_{k,1}+T_{k,2},
\end{align*}
where $T_{k,1}=k^{\tau/2} \Phi_k \bar {\mathbb{U}}_k \Big(\frac{\tilde {\mathbb{V}}_k-\bar {\mathbb{V}}_k}{\bar {\mathbb{V}}_k \tilde {\mathbb{V}}_k} \Big)+k^{\tau/2} \Phi_k \frac{\bar {\mathbb{U}}_k-\tilde {\mathbb{U}}_k}{\tilde {\mathbb{V}}_k}$ and $T_{k,2}=k^{\tau/2}\Phi_k\Big(E_{\theta_k}\Big[\frac{\tilde {\mathbb{U}}_k}{\tilde {\mathbb{V}}_k} \Big|\mathcal{F}_{k-1}\Big]-\frac{{\mathbb{U}}_k}{{\mathbb{V}}_k} \Big)$.
Note that
\begin{align}\label{tmp}\nonumber
\|T_{k,1}\|&\leq \|\Phi_k\| \frac{\|\bar {\mathbb{U}}_k \|}{|\bar {\mathbb{V}}_k \tilde {\mathbb{V}}_k|}k^{\tau/2}|\tilde {\mathbb{V}}_k-\bar {\mathbb{V}}_k |+\|\Phi_k\|\frac{1}{|\tilde {\mathbb{V}}_k|}k^{\tau/2}\|\bar {\mathbb{U}}_k-\tilde {\mathbb{U}}_k\|\\ \nonumber
&\leq \frac{\|\Phi_k\|\|\bar {\mathbb{U}}_k \|}{|\bar {\mathbb{V}}_k \tilde {\mathbb{V}}_k|}\frac{k^{\tau/2}}{N_k}\sum_{i=1}^{N_k}|\widehat S_{\theta_k}(H(x^i_k))-S_{\theta_k}(H(x^i_k)) |\\
&+\frac{\|\Phi_k\|}{|\tilde {\mathbb{V}}_k|}\frac{k^{\tau/2}}{N_k}\sum_{i=1}^{N_k}|\widehat S_{\theta_k}(H(x^i_k))-S_{\theta_k}(H(x^i_k)) |\|T(x^i_k)\|
\end{align}
Since $T(x)$ is bounded, it is easy to see that $\frac{\|\bar {\mathbb{U}}_k \|}{|\bar {\mathbb{V}}_k |}$ is also bounded. Furthermore, note that $|\tilde {\mathbb{V}}_k|$ is bounded away from zero. This, together with the boundedness of $\|\Phi_k\|$ and Assumption 3, imply that
the right-hand-side of (\ref{tmp}) converges to zero w.p.1.

For term $T_{k,2}$, let $\tilde {\mathbb{U}}_k^i$ and ${\mathbb{U}}_k^i$ be the $i$th components of $\tilde {\mathbb{U}}_k$ and ${\mathbb{U}}_k$, respectively. By using a second order two variable Taylor expansion of $\frac{\tilde {\mathbb{U}}_k^i}{\tilde {\mathbb{V}}_k}$ around $\frac{{\mathbb{U}}_k^i}{{\mathbb{V}}_k}$, we have
$$\frac{\tilde {\mathbb{U}}_k^i}{\tilde {\mathbb{V}}_k}=\frac{{\mathbb{U}}_k^i}{{\mathbb{V}}_k}+\frac{1}{{\mathbb{V}}_k}(\tilde {\mathbb{U}}_k^i-{\mathbb{U}}_k^i)-\frac{{\mathbb{U}}_k^i}{{\mathbb{V}}_k^2}(\tilde {\mathbb{V}}_k-{\mathbb{V}}_k)+\frac{\hat {\mathbb{U}}_k^i}{\hat {\mathbb{V}}_k^3}(\tilde {\mathbb{V}}_k-{\mathbb{V}}_k)^2-\frac{1}{\hat {\mathbb{V}}_k^2}(\tilde {\mathbb{U}}_k^i-{\mathbb{U}}_k^i)(\tilde {\mathbb{V}}_k-{\mathbb{V}}_k),$$
where $\hat {\mathbb{U}}_k^i$ and $\hat {\mathbb{V}}_k$ are on the line segments from $\tilde {\mathbb{U}}_k^i$ to ${\mathbb{U}}_k^i$ and from $\tilde {\mathbb{V}}_k$ to ${\mathbb{V}}_k$. Taking conditional expectations at both sides of the above equation, we have
\begin{align}\label{tmp1}\nonumber
\Big| E_{\theta_k}\Big[\frac{\tilde {\mathbb{U}}_k^i}{\tilde {\mathbb{V}}_k} \Big|\mathcal{F}_{k-1} \Big]-\frac{{\mathbb{U}}_k^i}{{\mathbb{V}}_k}  \Big|&\leq E_{\theta_k}\Big[\frac{|\hat {\mathbb{U}}_k^i|}{|\hat {\mathbb{V}}_k^3 |}   (\tilde {\mathbb{V}}_k-{\mathbb{V}}_k)^2\Big|\mathcal{F}_{k-1}\Big]+E_{\theta_k}\Big[\frac{1}{\hat {\mathbb{V}}_k^2}|(\tilde {\mathbb{U}}_k^i-{\mathbb{U}}_k^i)(\tilde {\mathbb{V}}_k-{\mathbb{V}}_k)| \Big|\mathcal{F}_{k-1}\Big]\\
&\leq \mathcal{C}_1E_{\theta_k}\Big[(\tilde {\mathbb{V}}_k-{\mathbb{V}}_k)^2\Big|\mathcal{F}_{k-1}\Big]+\mathcal{C}_2E_{\theta_k}\Big[|(\tilde {\mathbb{U}}_k^i-{\mathbb{U}}_k^i)(\tilde {\mathbb{V}}_k-{\mathbb{V}}_k)| \Big|\mathcal{F}_{k-1}\Big]
\end{align}
for constants $\mathcal{C}_1>0$ and $\mathcal{C}_2>0$. Thus, a straightforward calculation shows that the right-hand-side of (\ref{tmp1}) is $O(N_k^{-1})$. Consequently, we have $T_{k,2}\rightarrow 0$ w.p.1. as $k\rightarrow \infty$ by taking \ez{$N_k=N_0k^{\zeta}$} with $\zeta>\tau/2$. This shows $T_k\rightarrow 0$ w.p.1. as desired.
\endproof

\noindent \proof{\textbf{Lemma~\ref{lem2}.}}
$E_{\theta_k}[W_k|\mathcal{F}_{k-1}]=0$ follows directly from the definition
of $W_k$. Again, we
let $\tilde {\mathbb{U}}_k^i$ and ${\mathbb{U}}_k^i$ be the $i$th components of $\tilde {\mathbb{U}}_k$ and ${\mathbb{U}}_k$, let $T_i(x)$ be the $i$th component of the sufficient statistic $T(x)$, and define $\Sigma^k_{i,j}$ as the $(i,j)$th entry of the matrix $E_{\theta_k}[W_kW_k^T|\mathcal{F}_{k-1}]$. By using a first order two variable Taylor expansion of $\frac{\tilde {\mathbb{U}}_k^i}{\tilde {\mathbb{V}}_k}$ around $\frac{{\mathbb{U}}_k^i}{{\mathbb{V}}_k}$, we have
\begin{equation}\label{eqn:taylor1}
\frac{\tilde {\mathbb{U}}_k^i}{\tilde {\mathbb{V}}_k}=\frac{{\mathbb{U}}_k^i}{{\mathbb{V}}_k}+\frac{1}{{\mathbb{V}}_k}(\tilde {\mathbb{U}}_k^i-{\mathbb{U}}_k^i)-\frac{{\mathbb{U}}_k^i}{{\mathbb{V}}_k^2}(\tilde {\mathbb{V}}_k-{\mathbb{V}}_k)+\mathcal{R}_k,
\end{equation}
where $\mathcal{R}_k$ is a reminder term. Therefore, $\Sigma^k_{i,j}$ can be expressed as
\begin{align*}
\Sigma^k_{i,j}=&k^{\tau-\alpha}E_{\theta_k}\Big[\Big(\frac{\tilde {\mathbb{U}}^i_k}{\tilde {\mathbb{V}}_k}- E_{\theta_k}\Big[\frac{\tilde {\mathbb{U}}^i_k}{\tilde {\mathbb{V}}_k}\Big|\mathcal{F}_{k-1} \Big] \Big)\Big( \frac{\tilde {\mathbb{U}}^j_k}{\tilde {\mathbb{V}}_k}- E_{\theta_k}\Big[\frac{\tilde {\mathbb{U}}^j_k}{\tilde {\mathbb{V}}_k}\Big|\mathcal{F}_{k-1} \Big]\Big)  \Big|\mathcal{F}_{k-1}\Big]\\
=&k^{\tau-\alpha}\frac{1}{{\mathbb{V}}_k^2}E_{\theta_k}[(\tilde {\mathbb{U}}^i_k-{\mathbb{U}}^i_k)(\tilde {\mathbb{U}}^j_k-{\mathbb{U}}^j_k) |\mathcal{F}_{k-1}]~~~~[i]\\
&-k^{\tau-\alpha}\frac{{\mathbb{U}}^j_k}{{\mathbb{V}}_k^3}E_{\theta_k}[(\tilde {\mathbb{U}}^i_k-{\mathbb{U}}^i_k)(\tilde {\mathbb{V}}_k-{\mathbb{V}}_k) |\mathcal{F}_{k-1}]~~~~[ii]\\
&-k^{\tau-\alpha}\frac{{\mathbb{U}}^i_k}{{\mathbb{V}}_k^3}E_{\theta_k}[(\tilde {\mathbb{U}}^j_k-{\mathbb{U}}^j_k)(\tilde {\mathbb{V}}_k-{\mathbb{V}}_k) |\mathcal{F}_{k-1}]~~~~[iii]\\
&+k^{\tau-\alpha}\frac{{\mathbb{U}}^i_k{\mathbb{U}}^j_k}{{\mathbb{V}}_k^4}E_{\theta_k}[(\tilde {\mathbb{V}}_k-{\mathbb{V}}_k)^2 |\mathcal{F}_{k-1}]~~~~[iv]\\
&+k^{\tau-\alpha}\bar {\mathcal{R}}_k,
\end{align*}
where $\bar{\mathcal{R}}_k$ represents a higher-order term.
\begin{align*}
[i]&=k^{\tau-\alpha}\frac{1}{{\mathbb{V}}^2_k}\Big(E_{\theta_k}[\tilde {\mathbb{U}}^i_k \tilde {\mathbb{U}}^j_k |\mathcal{F}_{k-1}]-{\mathbb{U}}^i_k {\mathbb{U}}^j_k  \Big)\\
&=k^{\tau-\alpha}\frac{1}{{\mathbb{V}}^2_k}\frac{1}{N_k}\Big(E_{\theta_k}\big[S^2_{\theta_k}(H(X))T_i(X)T_j(X)  \big|\mathcal{F}_{k-1}\big]-  {\mathbb{U}}^i_k {\mathbb{U}}^j_k  \Big)\\
&=k^{\tau-\alpha}\frac{1}{N_k}\Big(\frac{E_{\theta_k}\big[S^2_{\theta_k}(H(X))T_i(X)T_j(X)  \big|\mathcal{F}_{k-1}\big]}{E_{\theta_k}^2[S_{\theta_k}(H(X))]}-\frac{{\mathbb{U}}^i_k{\mathbb{U}}^j_k}{{\mathbb{V}}^2_k}   \Big)\\
&=\frac{k^{\tau-\alpha}}{N_k}\Big[E_{p_k}\Big[T_i(X)T_j(X)\frac{p_k(X)}{f(X;\theta_k)} \Big]-E_{p_k}[T_i(X)]E_{p_k}[T_j(X)]  \Big].
\end{align*}
By using a similar argument, it can be seen that
\begin{align*}
[ii]&=\frac{k^{\tau-\alpha}}{N_k}\Big[E_{p_k}\big[T_j(X)\big]E_{p_k}\Big[T_i(X)\frac{p_k(X)}{f(X;\theta_k)} \Big]-E_{p_k}[T_i(X)]E_{p_k}[T_j(X)]  \Big],\\
[iii]&=\frac{k^{\tau-\alpha}}{N_k}\Big[E_{p_k}\big[T_i(X)\big]E_{p_k}\Big[T_j(X)\frac{p_k(X)}{f(X;\theta_k)} \Big]-E_{p_k}[T_i(X)]E_{p_k}[T_j(X)]  \Big],\\
[iv]&=\frac{k^{\tau-\alpha}}{N_k}\Big[E_{p_k}\big[T_j(X)\big]E_{p_k}\big[T_i(X)\big]E_{p_k}\Big[\frac{p_k(X)}{f(X;\theta_k)} \Big]-E_{p_k}[T_i(X)]E_{p_k}[T_j(X)]  \Big].
\end{align*}
Therefore,
\begin{align*}
\Sigma^k_{i,j}&=[i]-[ii]-[iii]+[iv]+k^{\tau-\alpha}\bar {\mathcal{R}}_k\\
&=\frac{k^{\tau-\alpha}}{N_k}E_{p_k}\Big[(T_i(X)-E_{p_k}[T_i(X)])(T_j(X)-E_{p_k}[T_j(X)])\frac{p_k(X)}{f(X;\theta_k)}  \Big]+k^{\tau-\alpha}\bar {\mathcal{R}}_k\\
&=\frac{k^{\tau-\alpha}}{N_k}E_{\theta_k}\Big[(T_i(X)-E_{p_k}[T_i(X)])(T_j(X)-E_{p_k}[T_j(X)])\frac{p_k^2(X)}{f^2(X;\theta_k)}  \Big]+k^{\tau-\alpha}\bar {\mathcal{R}}_k.
\end{align*}
By taking \ez{$N_k=N_0k^{\tau-\alpha}$}, it can be shown that the higher-order term
$k^{\tau-\alpha}\bar {\mathcal{R}}_k$ is $o(1)$. In addition, since $S_{\theta}(y)$ is continuous in $\theta$ for a fixed $y$, the point-wise convergence of $f(\cdot;\theta_k)$ to $f(\cdot;\theta^*)$ implies that $p_k(x)$ will also converge in a point-wise manner to a limiting distribution $p_*(x)$. Thus, the dominated convergence theorem suggests that
$\Sigma^k_{i,j}$ will converge to
$$\Sigma_{i,j}=\mathcal{C}E_{\theta^*}\Big[(T_i(X)-E_{p_*}[T_i(X)])(T_j(X)-E_{p_*}[T_j(X)])\frac{p_*^2(X)}{f^2(X;\theta^*)}  \Big]$$
for some positive constant $\mathcal{C}$. Therefore, the limiting matrix $\Sigma$ is given by
$$\Sigma=\mbox{Cov}_{\theta^*}\Big((T(X)-E_{p_*}[T(X)])\frac{p_*(X)}{f(X;\theta^*)}  \Big),$$
where $\mbox{Cov}_{\theta^*}(\cdot)$ is the covariance matrix with respect to $f(\cdot;\theta^*)$.

To show the last statement, we use H\"{o}lder's inequality and write
\begin{equation}\label{eqn:tmp}
\lim_{k\rightarrow \infty}E[I\{\|W_k \|^2\geq rk^{\alpha}\}\|W_k \|^2]
\leq \limsup_{k\rightarrow \infty}\Big[P\big(\|W_k\|^2\geq r k^{\alpha} \big) \Big]^{\frac{1}{2}}\Big[E\big[\|W_k \|^{4} \big]\Big]^{\frac{1}{2}}.\end{equation}
Note that
\begin{align*}
P\big(\|W_k\|^2\geq r k^{\alpha} \big)&=P\big(\|W_k\|\geq \sqrt{r} k^{\alpha/2} \big)\\
&\leq \frac{E[\|W_k\|^2]}{r k^{\alpha}}~~\mbox{by Chebyshev's inequality}\\
&=\frac{E\big[E_{\theta_k}[\|W_k\|^2|\mathcal{F}_{k-1}]\big]}{r k^{\alpha}}\\
&=\frac{E\big[\mbox{tr}(\Sigma^k)\big]}{r k^{\alpha}}\\
&=O(k^{-\alpha})
\end{align*}
by taking \ez{$N_k=N_0k^{\tau-\alpha}$} for $k$ sufficiently large, where the last step follows because all entries in $\Sigma^k$ are bounded and thus convergence w.p.1. implies convergence in expectation. On the other hand, by (\ref{eqn:taylor1}),
$E[\|W_k\|^{4}]$ can be expressed in terms of the fourth order central moments of the sample mean and it can be verified that $E[\|W_k\|^{4}]=O(1)$. This shows that the right-hand-side of (\ref{eqn:tmp}) is $O(k^{-\frac{\alpha}{2}})$, which vanishes to zero as $k\rightarrow \infty$.
\endproof

%
%


Acknowledgments: The authors gratefully acknowledge the support by the National Science Foundation under Grants ECCS-0901543 and CMMI-1130273 and Air Force Office of Scientific Research under YIP Grant FA-9550-12-1-0250. We are grateful to Xi Chen, graduate student in the Department of Industrial $\&$ Enterprise Systems Engineering at UIUC, for her help with conducting the numerical experiments in Section \ref{s:numerical}. 


\bibliographystyle{plain} 
\bibliography{Zhou-Bibtex}

\end{document}